# A Bayesian numerical homogenization method for elliptic multiscale inverse problems

Assyr Abdulle[*]    Andrea Di Blasio[†]

**Abstract**

A new strategy based on numerical homogenization and Bayesian techniques for solving multiscale inverse problems is introduced. We consider a class of elliptic problems which vary at a microscopic scale, and we aim at recovering the highly oscillatory tensor from measurements of the fine scale solution at the boundary, using a coarse model based on numerical homogenization and model order reduction. We provide a rigorous Bayesian formulation of the problem, taking into account different possibilities for the choice of the prior measure. We prove well-posedness of the effective posterior measure and, by means of G-convergence, we establish a link between the effective posterior and the fine scale model. Several numerical experiments illustrate the efficiency of the proposed scheme and confirm the theoretical findings.

## 1 Introduction

Inverse problems for partial differential equations (PDEs) play an important role in the sciences and the engineering, with numerous applications as geoscience or medical imaging for example. In this work we are interested in PDEs characterized by the presence of variations on a very fine scale, which can be found for example in the study of composite materials or pourous media. Let $\Omega \in \mathbb{R}^d$, $d \geq 2$, be an open, bounded, connected set with sufficiently regular boundary $\partial \Omega$, and consider the problem of finding $u^\varepsilon \in H^1(\Omega)$ such that

$$\begin{aligned} -\nabla \cdot (A^\varepsilon \nabla u^\varepsilon) &= 0 &&\text{in } \Omega, \\ u^\varepsilon &= g &&\text{on } \partial \Omega. \end{aligned} \qquad (1)$$

The tensor $A^\varepsilon = A^\varepsilon(x)$, $x \in \Omega$, belongs to $M(\alpha, \beta, \Omega)$, $0 < \alpha < \beta$, where

$M(\alpha, \beta, \Omega) :=$
$\{A \in L^\infty(\Omega, Sym_d) \,:\, \alpha |\xi|^2 \leq A(x)\xi \cdot \xi, \, |A(x)\xi| \leq \beta |\xi| \, \forall \xi \in \mathbb{R}^d \text{ and a. e. in } \Omega\},$

and $Sym_d$ denotes the class of $d \times d$ symmetric real valued matrices. The superscript in $A^\varepsilon$ (respectively $u^\varepsilon$) emphasizes that the tensor (the solution) varies on a fine scale proportional to $\varepsilon$, which is usually much smaller than the domain $\Omega$ considered for application. The inverse problem we are interested in, is to recover the highly oscillatory tensor $A^\varepsilon$ based on observations originating from (1). Often standard numerical techniques such as the Finite Element Method (FEMs) are not appropriate to approximate (1) since mesh resolution at the finest scale is required to provide a reliable solution. Mesh resolution down to the $\varepsilon$ scale can be prohibitively expensive when $\varepsilon$ is small. This issue is even more dramatic when solving inverse problems, since one typically needs multiple evaluations of (1), and thus an alternative approach is required. From the theory of homogenization [11,13,17] we know that there exists an effective tensor $A^0$ such that (up to a subsequence) the solution of (1) converges in a weak sense to the solution $u^0 \in H^1(\Omega)$ of the problem

$$
\begin{aligned}
-\nabla \cdot (A^0 \nabla u^0) &= 0 \quad \text{in } \Omega\,, \\
u^0 &= g \quad \text{on } \partial\Omega\,,
\end{aligned}
\tag{2}
$$

where $A^0$ is referred to as the homogenized tensor. An explicit form of $A^0$ is usually not known, and so numerical homogenization [1,2,5] is needed to obtain the homogenized solution $u^0$ based on data defining problem (1). Our strategy to efficiently retrieve the conductivity $A^\varepsilon$, based on observations originating from (1), relies on the reduced model (2). In [4] we analyzed and solved the inverse problem in the case where the observed quantities were defined by the Dirichlet to Neumann map associated to (1),

$$\Lambda_{A^\varepsilon} : g \in H^{1/2}(\partial\Omega) \mapsto A^\varepsilon \nabla u^\varepsilon \cdot \nu|_{\partial\Omega} \in H^{-1/2}(\partial\Omega)\,, \tag{3}$$

where $\nu$ denotes the exterior unit normal to $\partial\Omega$.

---

[*]École Polytechnique Fédérale de Lausanne, Switzerland (assyr.abdulle@epfl.ch)
[†]École Polytechnique Fédérale de Lausanne, Switzerland(andrea.diblasio@epfl.ch)

In this paper, as in [4], we consider a class of parametrized multiscale locally periodic tensors of the form $A^\varepsilon_{\sigma^*}(x) = A(\sigma^*(x), x/\varepsilon)$, where $\sigma^* : \Omega \to \mathbb{R}$. We assume that the map $(t,x) \mapsto A(t, x/\varepsilon)$, $t \in [\sigma^-, \sigma^+]$, $x \in \Omega$, is known while $\sigma^* : \Omega \to \mathbb{R}$ has to be determined to recover the whole tensor. A typical example of this setting could be represented by a multi-phase medium, whose constituent materials are known, but their respective volume fraction or marcoscopic orientation are unknown. Departing from [4], where in order to ensure well-posedness we solved the problem by means of Tikhonov regularization, we recast here the problem into a statistical framework, and develop a multiscale numerical method based on Bayesian techniques. In addition in contrast to [4], instead of considering observed data as living in some functional space, e.g. $H^{-1/2}(\partial\Omega)$, we consider discrete quantities in $\mathbb{R}$ represented by the average of the normal flux at the boundary measured on different locations $\Gamma_j \subset \partial\Omega$, $j = 1, \ldots, J$, $J \in \mathbb{N}$. For a survey on the Bayesian approach for inverse problems, we mention [12,19]. For a rigorous Bayesian formulation of the inverse conductivity problem, known also as electrical impedance tomography (EIT), we also mention [14]. We mention that Bayesian multiscale inverse problems have also been addressed in [18]. The contribution of this paper can be summarized as follows:

- because of the prohibitive cost of the forward problem in a multiscale context, we introduce an *effective forward problem* and a related *effective posterior measure*. The modeling error introduced in this framework can be quantified in terms of $G-convergence$ and we provide an offline algorithm to correct for the model discrepancy;

- our numerical algorithm makes use of multiscale methods and model order reduction techniques to tackle computationally challenging multi-dimensional multiscale problems;

- our methodology allows to effectively recover a multiscale conductivity tensor through partial observations on the boundary of the domain.

Following [12,19], we give a rigorous Bayesian formulation of the multiscale problem, and prove the well-posedness of the effective posterior measure for our setting. We employ different kind of prior measures, considering log-Gaussian and level set priors. Moreover, we establish a link between the effective posterior measure and the fine scale model in terms of Hellinger distance, using G-convergence, to quantify the discrepancy between the homogenized data and the data originating from (1). The numerical method builds on the reduced basis heterogenenous multiscale method developed in [3]. Finally, inspired by [9], we approximate numerically the modelling error distribution and we verify that including the modelling error distribution in the definition of the posterior measure can improve significantly the results, especially when $\varepsilon$ is relatively large.

The outline of the work is as follows. In Section 2 we describe our setting for the observed data and we recall some useful tools for the Bayesian approach to inverse problems. In Section 3 we state some preliminary results on well-posedness of the posterior measure and homogenization theory and we introduce two types of prior measures that will be used throughout the work. Our main results are presented in Section 4. We prove existence and well-posedness of the effective posterior, and establish the convergence of the Hellinger distance between the effective posterior and the posterior measure based on the full fine scale model. In Section 5 we give a brief survey on the Markov chain Monte Carlo method used to sample from the posterior distribution, while in Section 6 we explain how to approximate numerically (2) by a model order reduction multiscale method. Numerical experiments that illustrate our multiscale inverse method and confirm our theoretical findings are presented in Section 7.

## 2 Preliminaries: problem definition, homogenization and G-convergence

Let $\Omega$ be an open and bounded set in $\mathbb{R}^d$. We consider a class of parametrized multiscale locally periodic tensors of the type $A^\varepsilon_{\sigma^*}(x) = A(\sigma^*(x), x/\varepsilon)$, where $\sigma^* : \Omega \to \mathbb{R}$. Our aim is to recover $A^\varepsilon_{\sigma^*}$ from measurements originating from the model

$$\begin{aligned} -\nabla \cdot (A^\varepsilon_{\sigma^*} \cdot \nabla u^\varepsilon) &= 0 &&\text{in } \Omega, \\ u^0 &= g &&\text{on } \partial\Omega. \end{aligned} \qquad (4)$$

Our unknown is represented by $\sigma^*$, while we assume to know the map $(t, x) \mapsto A(t, x/\varepsilon)$, $t \in [\sigma^-, \sigma^+]$, $x \in \Omega$. Let us consider then the following admissible set for the unknown $\sigma^* : \Omega \to \mathbb{R}$

$$U := \{\sigma \in L^\infty(\Omega) \,:\, \sigma^- \leq \sigma(x) \leq \sigma^+\}.$$

We consider $J \in \mathbb{N}$ boundary portions of $\partial\Omega$, and we denote them as $\Gamma_j \in \partial\Omega$, $j = 1, \ldots, J$, $\Gamma_i \cap \Gamma_j = \emptyset$ for $i \neq j$. These portions of the boundary represents the locations at which the measurements are carried out. Moreover the same experiment is reproduced for $L \in \mathbb{N}$ different Dirichlet data, which we denote as $g_l$, $l = 1, \ldots, L$. Hence we have $J \times L$ observations. Then we may introduce the forward operator $F^\varepsilon : U \to \mathbb{R}^{JL}$, $F^\varepsilon(\sigma) = \text{vec}(\{f^\varepsilon_{jl}(\sigma)\}_{\substack{1 \leq j \leq J \\ 1 \leq l \leq L}})$,

$$f^\varepsilon_{jl}(\sigma) = \int_{\Gamma_j} \Lambda_{A^\varepsilon_\sigma} g_l \cdot \phi_j \, ds, \quad j = 1, \ldots, J, \ l = 1, \ldots, L, \qquad (5)$$

where $\Lambda_{A^\varepsilon_\sigma}$ is the Dirichlet to Neumann map (3) associated to the tensor $A^\varepsilon_\sigma(x) = A(\sigma(x), x/\varepsilon)$, and $\phi_j \in H^{1/2}(\partial\Omega)$ such that $supp(\phi_j) \subseteq \Gamma_j$ for all $j = 1, \ldots, J$. In the following setting, we assume to

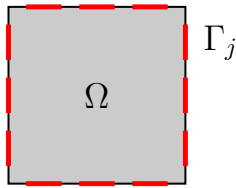

Figure 1: Picture representing the computational domain $\Omega$ and the boundary portions $\Gamma_j$ used to compute the observations.

dispose of a finite number of observations, corrupted by some noise, so that

$$z = F^\varepsilon(\sigma^*) + e, \quad e \sim \mathcal{N}(0, C_e), \qquad (6)$$

where $C_e$ is a given covariance matrix. Based on these measurements we would like to recover $\sigma^*$. Let $X$ be a Banach space, and $P$ some map $P : \theta \in X \mapsto \sigma \in U$. The introduction of $X$ and $P$ will be useful later on to build different kind of prior measures on the admissible set $U$. Introducing this abstract framework is also useful to perform a rigorous analysis about the validity of our approach, which will be carried out in Section 3. Let us define the potential function $\Phi^\varepsilon : X \times \mathbb{R}^{JL} \to \mathbb{R}$, which measures the distance between the observed data and the values produced by the observation model for some $\theta \in X$ as

$$\Phi^\varepsilon(\theta, z) = \frac{1}{2}||z - G^\varepsilon(\theta)||^2_{C_e} \tag{7}$$
$$= \frac{1}{2}(z - G^\varepsilon(\theta))^\top C_e^{-1}(z - G^\varepsilon(\theta)),$$

where $G^\varepsilon = F^\varepsilon \circ P$. Simply trying to minimize (7) leads to an ill posed problem. To ensure well posedness we may add some regularization term (e.g. Tikhonov regularization) or recast the problem into a statistical framework, where all the quantities involved are treated as random variables (Bayesian approach). Differently from standard regularization techniques which produce as solution a single point estimate of the unknown, with the statistical approach the solution is represented by a probability measure, so called the posterior probability measure. The posterior measure can then be used to infer about the parameter values and quantify their uncertainties. In Bayesian theory, it is assumed that all the prior informations we dispose about the unknown we are seeking for, can be described by what is called the prior measure, which we denote here as $\mu_{pr}$. Using (6) and applying Bayes' formula we obtain that the posterior measure of $\theta$ given $z$, denoted as $\mu^\varepsilon(\theta|z)$, is related to $\mu_{pr}$ through the Radon-Nikodym derivative

$$\frac{\mathrm{d}\mu^\varepsilon(\theta|z)}{\mathrm{d}\mu_{pr}(\theta)} \propto \exp(-\Phi^\varepsilon(\theta, z)). \tag{8}$$

Unfortunately trying to explore $\mu^\varepsilon(\theta|z)$ via sampling techniques as Markov chain Monte Carlo methods (MCMC) is infeasible, due to the high computational effort needed to evaluate the model $G^\varepsilon$ even for few realizations of $\theta \in X$. Hence, to drastically reduce the computational cost we combine the inverse problem with a coarse graining strategy. To do so let us recall briefly some results from homogenization theory [7, 13, 17], in particular the concept of G-convergence.

**Definition 1.** *Let $\{A^\varepsilon\}_{\varepsilon>0}$ be a sequence of matrices in $M(\alpha, \beta, \Omega)$. We say that $\{A^\varepsilon\}_{\varepsilon>0}$ G-converges to the matrix $A^0 \in M(\alpha, \beta, \Omega)$ if and only if for every function $f \in H^{-1}(\Omega)$, $g \in H^{1/2}(\partial\Omega)$, the solution $u^\varepsilon$ of*

$$\begin{aligned} -\nabla \cdot (A^\varepsilon \nabla u^\varepsilon) &= f &&\text{in } \Omega, \\ u^\varepsilon &= g &&\text{on } \partial\Omega. \end{aligned} \tag{9}$$

*is such that*
$$u^\varepsilon \rightharpoonup u^0 \quad \text{weakly in } H^1(\Omega),$$

*where $u^0$ is the unique solution of*

$$\begin{aligned} -\nabla \cdot (A^0 \nabla u^0) &= f &&\text{in } \Omega, \\ u^0 &= g &&\text{on } \partial\Omega. \end{aligned} \tag{10}$$

*A consequence of G-convergence is the weak convergence of the flux*

$$A^\varepsilon \nabla u^\varepsilon \rightharpoonup A^0 \nabla u^0 \quad \text{weakly in } (L^2(\Omega))^d.$$

**Theorem 1.** *(See for example [11, 17]). One has the following compactness result. Let $\{A^\varepsilon\}_{\varepsilon>0}$ be a sequence of matrices in $M(\alpha, \beta, \Omega)$. Then there exists a subsequence $\{A^{\varepsilon'}\}_{\varepsilon'>0}$ and a matrix $A^0 \in M(\alpha, \beta, \Omega)$ such that $\{A^{\varepsilon'}\}_{\varepsilon'>0}$ G-converges to $A^0$.*

In particular we consider for a $\sigma \in U$ the sequence of $Y$-periodic matrices defined by

$$A^\varepsilon_\sigma(x) = A(\sigma(x), x/\varepsilon) = A(\sigma(x), y), \quad A(\sigma(x), \cdot) \in M(\alpha, \beta, Y), \forall x \in \Omega,$$

$$A^\varepsilon_\sigma(x) = \{a^\varepsilon_{ij}(x)\}_{1 \le i, j \le d} \quad \text{a.e. on } \mathbb{R}^d,$$

where

$$a_{ij}^{\varepsilon}(x) = a_{ij}(\sigma(x), x/\varepsilon) = a_{ij}(\sigma(x), y), a_{ij}(\sigma(x), \cdot) \text{ is } Y\text{-periodic}, \forall x \in \Omega, \forall i, j = 1, \ldots, d,$$

and $Y$ denotes the reference unit cell $(0,1)^d$. Such tensors are usually referred to as locally periodic in the literature. In this particular case we have that the whole sequence $\{A_\sigma^\varepsilon\}_{\varepsilon>0}$ G-converges to the tensor $A_\sigma^0 \in M(\alpha, \beta, \Omega)$, $A_\sigma^0(x) = A^0(\sigma(x)) = \{a_{ij}^0(\sigma(x))\}_{1 \leq i,j \leq d}$, which is elliptic and is given by

$$a_{ij}^0(\sigma(x)) = \frac{1}{|Y|} \int_Y a_{ij}(\sigma(x), y) \, \mathrm{d}y - \frac{1}{|Y|} \sum_{k=1}^d \int_Y a_{ik}(\sigma(x), y) \frac{\partial \chi_j}{\partial y_k} \, \mathrm{d}y \quad \forall \, i, j = 1, \ldots, d.$$

The micro functions $\chi_j$, $j = 1, \ldots, d$, are defined to be the unique solutions of the cell problems: find $\chi_j \in W_{per}^1(Y)$ such that

$$\int_Y A(\sigma(x), y) \nabla_y \chi_j \cdot \nabla_y v \, \mathrm{d}y = \int_Y A(\sigma(x), y) \mathbf{e_j} \cdot \nabla_y v \, \mathrm{d}y, \forall v \in W_{per}^1(Y), \tag{11}$$

where $\{\mathbf{e_j}\}_{j=1}^d$ is the canonical basis of $\mathbb{R}^d$ and

$$W_{per}^1(Y) = \{v \in H_{per}^1(Y) : \int_Y v \, \mathrm{d}y = 0\},$$

and $H_{per}^1(Y)$ is defined as the closure of $C_{per}^\infty(Y)$ for the $H^1$-norm (where $C_{per}^\infty(Y)$ denotes the subset of $C^\infty(\mathbb{R}^d)$ of periodic functions in $Y$).

Hence using homogenization theory, we may introduce the operator $F^0 : U \to \mathbb{R}^{JL}$, $F^0(\sigma) = \text{vec}(\{f_{jl}^0(\sigma)\}_{\substack{1 \leq j \leq J \\ 1 \leq l \leq L}})$,

$$f_{jl}^0(\sigma) = \int_{\Gamma_j} \Lambda_{A_\sigma^0} g_l \cdot \phi_j \, \mathrm{d}s, \; j = 1, \ldots, J, \; l = 1, \ldots, L, \tag{12}$$

where $\Lambda_{A_\sigma^0}$ is the Dirichlet to Neumann map associated to the tensor $A_\sigma^0$, the homogenized tensor associated to $A_\sigma^\varepsilon$. Then we can define a new potential function $\Phi^0 : X \times \mathbb{R}^{JL} \to \mathbb{R}$, as

$$\Phi^0(\theta, z) = \frac{1}{2} \|z - G^0(\theta)\|_{C_e}^2, \tag{13}$$

$G^0 : F^0 \circ P$, where $P$ is a map such that $P : X \to U$. As for the full fine scale model we can invoke Bayes' formula to define a posterior measure $\mu^0(\theta|z)$ associated to the potential function (13) which satisfies

$$\frac{\mathrm{d}\mu^0(\theta|z)}{\mathrm{d}\mu_{pr}(\theta)} \propto \exp(-\Phi^0(\theta, z)). \tag{14}$$

We note that this new measure is much easier to explore via sampling techniques since the homogenized forward model $F^0 : U \to \mathbb{R}^{JL}$ can be approximated efficiently and indipendently on $\varepsilon$.

## 3  Well-posedness of the posterior measure

In this section we recall some theoretical results about existence and well-posedness of the posterior measure. It is important to underline that existence and well-posedness of the posterior measure is typically determined from continuity properties of the forward operator entering in the definition of the potential function. Then it is necessary to build prior measures such that every proposal lies in the function space on which the continuity properties of the forward operator are satisfied. Hence, some analysis on regularity properties of the forward operator is needed. This is carried on in what follows. We assume to have a prior Gaussian measure $\mu_{pr} = \mathcal{N}(\theta_{pr}, C_{pr})$ defined on a Banach space $X$. Let $\mu^0(\theta|z)$ be a posterior measure that we assume as in Section 2 to satisfy

$$\frac{\mathrm{d}\mu^0(\theta|z)}{\mathrm{d}\mu_{pr}(\theta)} = \frac{1}{C^0(z)} \exp(-\Phi^0(\theta, z)), \tag{15}$$

where $\Phi^0(\theta, z)$ is the potential defined in (13) and $C^0(z)$ is the normalization constant

$$C^0(z) = \int_X \exp(-\Phi^0(\theta, z))\mu_{pr}(\mathrm{d}\theta),$$

so that $\mu^0(\theta|z)$ is actually a probability measure.

**Definition 2.** *Let $\mu^1$ and $\mu^2$ be two probability measures on a Banach space $X$. Assume $\mu^1$ and $\mu^2$ are both absolutely continuous with respect to a common reference measure $\mu$, defined on the same measure space. Then the Hellinger distance between $\mu^1$ and $\mu^2$ is defined as*

$$d^2_{Hell}(\mu^1, \mu^2) = \frac{1}{2}\int_X \left(\sqrt{\frac{\mathrm{d}\mu^1}{\mathrm{d}\mu}} - \sqrt{\frac{\mathrm{d}\mu^2}{\mathrm{d}\mu}}\right)^2 \mathrm{d}\mu.$$

The next theorem gives sufficient conditions on $\Phi^0 : X \times \mathbb{R}^{JL} \to \mathbb{R}$ and $\mu_{pr}$ for the posterior measure defined in (15) to be well-defined.

**Theorem 2.** *(Theorem 2.3 from [14]). Assume the function $\Phi^0 : X \times \mathbb{R}^{JL} \to \mathbb{R}$ and the probability measure $\mu_{pr}$ on the probability space $(X, \Sigma)$ satisfy the following properties:*

1. *For every $r > 0$ there is a $K = K(r)$ such that for all $\theta \in X$ and for all $z \in \mathbb{R}^{JL}$ such that $\max\{||\theta||_X, ||z||_{C_e}\} < r$*

$$0 \leq \Phi^0(\theta, z) \leq K.$$

2. *For any fixed $z \in \mathbb{R}^{JL}$ the function $\Phi^0(\cdot, z) : X \to \mathbb{R}$ is continuous $\mu_{pr}$-almost surely on the probability space $(X, \Sigma, \mu_{pr})$.*

3. *For $z_1, z_2 \in \mathbb{R}^{JL}$ with $\max\{||z_1||_{C_e}, ||z_2||_{C_e}\} < r$ and for every $\theta \in X$, there exists $C = C(r, ||\theta||_X)$ such that*

$$|\Phi^0(\theta, z_1) - \Phi^0(\theta, z_2)| \leq C||z_1 - z_2||_{C_e}.$$

*Then the posterior measure $\mu^0(\theta|z)$ given by (15) is a well-defined probability measure and it is Lipschitz in the data $z$, with respect to the Hellinger distance: if $\mu^0(\theta|z_1)$ and $\mu^0(\theta|z_2)$ are two measures corresponding to data $z_1$ and $z_2$, then there is a constant $C = C(r) > 0$ such that, for all $z_1, z_2$ with $\max\{||z_1||_{C_e}, ||z||_{C_e}\} < r$,*

$$d_{Hell}(\mu^0(\theta|z_1), \mu^0(\theta|z_2)) \leq C||z_1 - z_2||_{C_e}.$$

We consider the case where $\mu_{pr}$ is a Gaussian probability measure on the Banach space $X = C^0(\overline{\Omega})$, and we will show in Section 4 that assumptions 1-3 are satisfied by $\mu_{pr}$ and $G^0 = F^0 \circ P$, where $P : \theta \in C^0(\overline{\Omega}) \mapsto \sigma \in U$ is some map such that if $||\theta - \theta_n||_{L^\infty(\Omega)} \to 0$, then $P(\theta_n) \to P(\theta)$ either uniformly or in measure. In particular we consider two different definitions of $P$, which we denote as $P_1$ and $P_2$, described in what follows.

## 3.1 Prior measure

The map $P_1$ is simply defined as $P_1(\theta) = \exp(\theta)$. By continuity of $P_1$ we see that if $\theta \in C^0(\overline{\Omega})$ and $\{\theta_n\}_{n>0} \subset C^0(\overline{\Omega})$ is a sequence such that $||\theta - \theta_n||_{L^\infty(\Omega)} \to 0$, then $||P_1(\theta) - P_1(\theta_n)||_{L^\infty(\Omega)} \to 0$. Moreover note that since $\theta$ is distributed according to a Gaussian measure, $P_1(\theta)$ is distributed according to a log-Gaussian measure.

The map $P_2$, which in [16] is referred to as level set prior, is defined instead in the following way. Let $n \in \mathbb{N}$ and fix constants $-\infty = c_0 < \ldots < c_n = \infty$. Given $\theta : \Omega \to \mathbb{R}$, define $\Omega_i \subseteq \Omega$ as

$$\Omega_i = \{x \in \Omega : c_{i-1} \leq \theta(x) < c_i\}, \quad i = 1, \ldots, n,$$

so that $\overline{\Omega} = \cup_{i=1}^n \overline{\Omega}_i$ and $\Omega_i \cap \Omega_j = \emptyset$ for $i \neq j$. Let us also define the level sets

$$\Omega_i^0 = \overline{\Omega}_i \cap \overline{\Omega}_{i+1} = \{x \in \Omega : \theta(x) = c_i\}, \quad i = 1 \ldots, n-1.$$

Now given some strictly positive functions $f_1, \ldots, f_n \in C^0(\overline{\Omega})$, we define the map $P_2 : C^0(\overline{\Omega}) \to U$ as

$$P_2(\theta) = \sum_{i=1}^n f_i \mathbb{1}_{\Omega_i}.$$

In particular we will consider $f_i$ which are constant on $\Omega$. For the continuity of the map $P_2$, we have the following theorem (we denote by $|\Omega_j|$ the measure of $\Omega_j$).

**Proposition 1.** *(Proposition 2.6 and 2.8 in [16]). Let $\{\theta_n\}_{n>0} \subset C^0(\overline{\Omega})$ converge to some $\theta \in C^0(\overline{\Omega})$ uniformly. Then $\{P_2(\theta_n)\}_{n>0}$ converges to $P_2(\theta)$ in $L^q(\Omega)$, $1 \leq q < \infty$, if and only if $|\Omega_i^0| = 0$ for all $i = 1, \ldots, n-1$. Let $\mu_{pr}$ be a Gaussian probability measure on $C^0(\overline{\Omega})$ and let $\theta \sim \mu_{pr}$. Then $|\Omega_i^0| = 0$ $\mu_{pr}$-almost surely for $i = 1, \ldots, n-1$.*

## 4 Main results

Here we discuss our main contributions. We recall first a regularity result for the class of $d \times d$ symmetric matrix functions $t \to A^0(t)$ which was obtained in [4] based on continuity assumptions on the class of $d \times d$ (fine scale) symmetric matrix functions $(t, y) \to A(t, y)$, $y = x/\varepsilon$. Afterwards, we study the continuity of the forward operator $G^0 : X \to \mathbb{R}^{JL}$, and we analyse the modelling error switching from multiscale observations to the ones produced by the homogenized model.

**Theorem 3.** *(See [4]). Let $x/\varepsilon = y \in Y$. Consider the class of $d \times d$ symmetric matrix functions $(t, y) \mapsto A(t, y)$, where $a_{ij}$ is $Y$-periodic, $\forall i, j = 1, \ldots, d$, $t \in [\sigma^-, \sigma^+]$, $0 < \sigma^- < \sigma^+$. Assume that there exist $E_1$, $\alpha$, $\beta > 0$ such that*

$$A \in W^{1,\infty}([\sigma^-, \sigma^+] \times Y, Sym_d), \; ||A||_{W^{1,\infty}([\sigma^-,\sigma^+];W^{1,\infty}(Y))} \leq E_1. \tag{16}$$

$$\partial_t A \in W^{1,\infty}([\sigma^-, \sigma^+] \times Y, Sym_d), \; ||\partial_t A||_{W^{1,\infty}([\sigma^-,\sigma^+];W^{1,\infty}(Y))} \leq E_1. \tag{17}$$

$$\alpha|\xi|^2 \leq A(t,y)\xi \cdot \xi, |A(t,y)\xi| \leq \beta|\xi|, \quad \text{for a.e. } y \in Y \text{ and } \forall t \in [\sigma^-, \sigma^+], \xi \in \mathbb{R}^d. \tag{18}$$

$$\partial_t A(t,y)\xi \cdot \xi \geq E_1^{-1}|\xi|^2, \quad \text{for a.e. } y \in Y \text{ and } \forall t \in [\sigma^-, \sigma^+], \xi \in \mathbb{R}^d. \tag{19}$$

*Then there exist $E_2$, $\alpha$, $\beta > 0$ such that homogenized map $t \mapsto A^0(t)$ satisfies*

$$|\partial_t A^0(t)| + |\partial_t^2 A^0(t)| \leq E_2, \qquad \forall t \in [\sigma^-, \sigma^+]. \tag{20}$$

$$\alpha|\xi|^2 \leq A^0(t)\xi \cdot \xi, |A^0(t)\xi| \leq \beta|\xi|, \qquad \forall t \in [\sigma^-, \sigma^+], \xi \in \mathbb{R}^d. \tag{21}$$

$$\partial_t A^0(t)\xi \cdot \xi \geq E_2^{-1}|\xi|^2, \qquad \forall t \in [\sigma^-, \sigma^+], \xi \in \mathbb{R}^d. \tag{22}$$

### 4.1 Continuity of the forward operator

The following theorem establish the continuity of our forward operator.

**Theorem 4.** *Let the assumptions of Theorem 3 be satisfied. Let the sequence $\{\sigma_n\}_{n>0} \subset U$ converge to some $\sigma \in U$ either uniformly or in measure. Then the sequence $\{\Lambda_{A^0_{\sigma_n}} g\}_{n>0}$ converges to $\Lambda_{A^0_\sigma} g$ in $H^{-1/2}(\partial\Omega)$.*

*Proof.* The first part pf the result has been proved in [4]. For convenience we briefly recall the arguments. Let us define $w = A^0_\sigma \nabla u^0(\sigma) - A^0_{\sigma_n} \nabla(\sigma_n)$. Observing that $w \in H(\Omega, \text{div})$ and using the continuity of the map $w \in H(\Omega, \text{div}) \mapsto w \cdot \nu \in H^{-1/2}(\partial\Omega)$ we know that

$$||w \cdot \nu||_{H^{-1/2}(\partial\Omega)} \leq ||w||_{L^2(\Omega)}.$$

Using Cauchy-Schwarz inequality and (20)-(21) we obtain

$$\int_\Omega |w|^2 \, dx = \int_\Omega A^0_{\sigma_n}(\nabla u^0(\sigma) - \nabla u^0(\sigma_n)) \cdot w \, dx$$

$$+ \int_\Omega (A^0_\sigma - A^0_{\sigma_n}) \nabla u^0(\sigma) \cdot w \, dx$$

$$\leq \beta ||\nabla u^0(\sigma) - \nabla u^0(\sigma_n)||_{L^2(\Omega)} ||w||_{L^2(\Omega)}$$

$$+ C_2 \left( \int_\Omega |\sigma - \sigma_n|^2 |\nabla u^0(\sigma)|^2 \, dx \right)^{1/2} ||w||_{L^2(\Omega)}. \tag{23}$$

It follows from the weak formulation of $u^0(\sigma)$ and $u^0(\sigma_n)$ that, for all $v \in H^1_0(\Omega)$, we have that

$$\int_\Omega (A^0_\sigma \nabla u^0(\sigma) - A^0_{\sigma_n} \nabla u^0(\sigma_n)) \cdot \nabla v \, dx = 0.$$

Then

$$\int_\Omega A^0_\sigma (\nabla u^0(\sigma) - \nabla u^0(\sigma_n)) \cdot \nabla v \, dx = \int_\Omega (A^0_{\sigma_n} - A^0_\sigma) \nabla u^0(\sigma_n) \cdot \nabla v \, dx.$$

By choosing $v = u^0(\sigma) - u^0(\sigma_n) \in H^1_0(\Omega)$ and using Cauchy-Schwarz inequality, (20), (21) we obtain

$$||\nabla u^0(\sigma) - \nabla u^0(\sigma_n)||_{L^2(\Omega)} \leq \alpha^{-1} E_2 \left( \int_\Omega |\sigma - \sigma_n|^2 |\nabla u^0(\sigma_n)|^2 \, dx \right)^{1/2}. \tag{24}$$

Inserting (24) into (23) we obtain

$$||w||_{L^2(\Omega)} \leq E_2(1 + \alpha^{-1}\beta) \left( \int_\Omega |\sigma - \sigma_n|^2 |\nabla u^0(\sigma)|^2 \, dx \right)^{1/2}, \tag{25}$$

and by using Holder's inequality and Lax-Milgram we finally obtain

$$||w||_{L^2(\Omega)} \leq E_2(1 + \alpha^{-1}\beta) ||\sigma - \sigma_n||_{L^\infty(\Omega)} ||\nabla u^0(\sigma)||_{L^2(\Omega)}$$

$$\leq E_2 \alpha^{-1} \beta (1 + \alpha^{-1}\beta) ||g||_{H^{1/2}(\partial\Omega)} ||\sigma - \sigma_n||_{L^\infty(\Omega)}. \tag{26}$$

Now, if $||\sigma - \sigma_n||_{L^\infty(\Omega)} \to 0$ the result follows from (26). On the other hand if $|\sigma - \sigma_n| \to 0$ in measure, since $|\Omega| < \infty$ and $\nabla u^0(\sigma) \in (L^2(\Omega))^d$, we have also that the integrand of (25) $|\sigma - \sigma_n|^2 |\nabla u^0(\sigma)|^2 \to 0$ in measure (see corollary 2.2.6 in [8] for example). Now, since $|\sigma - \sigma_n|$ is uniformly bounded by assumptions, the whole integrand is bounded by a scalar multiple of $|\nabla u^0(\sigma)|^2$. Therefore by applying the Lebesgue's dominated convergence theorem, we obtain that $|\sigma - \sigma_n|^2 |\nabla u^0(\sigma)|^2 \to 0$ in $L^1(\Omega)$, and the result follows. □

*Remark 1.* The Lebesgue's dominated convergence theorem is stated for sequences converging almost everywhere. However, convergence almost everywhere can be replaced in this case by convergence in measure, since $|\Omega| < \infty$.

*Remark 2.* The result given in Theorem 4 is stronger than the one we obtained in Lemma 4.1 in [4], since it states continuity of the flux with respect to the convergence in measure of $\{\sigma_n\}_{n>0}$ to $\sigma$. There continuity of the flux with respect to the $L^r(\Omega)$ topology on $U$, with $1 \leq r < \infty$, was obtained by asking for higher regularity of the solution $u^0$.

We can deduce the following corollary that establishes the contnuity of the effective forward operator.

**Corollary 1.** *Let the assumptions of Theorem 3 be satisfied. Let the sequence $\{\sigma_n\}_{n>0} \subset U$ converge to some $\sigma \in U$ either uniformly or in measure. Then the sequence $\{F^0(\sigma_n)\}_{n>0}$ converges to $F^0(\sigma)$.*

*Proof.* We have that

$$||F^0(\sigma) - F^0(\sigma_n)||_{C_e} \leq C \sum_{j=1}^{J} \sum_{l=1}^{L} \int_{\Gamma_j} \left|(\Lambda_{A^0_\sigma} - \Lambda_{A^0(\sigma_n)})g_l \cdot \phi_j\right| \, \mathrm{d}s$$

$$\leq C \sup_l ||(\Lambda_{A^0_\sigma} - \Lambda_{A^0_{\sigma_n}})g_l||_{H^{-1/2}(\partial\Omega)} \sup_j ||\phi_j||_{H^{1/2}(\partial\Omega)},$$

and the result follows from Theorem 4. □

Finally we establish that the posterior measure (14) based on the potential function $\Phi^0$ is well defined and Lipschitz contnuous in the data with respect to the Hellinger distance.

**Corollary 2.** *Let the assumptions of Theorem 3 be satisfied. Let $\mu_{pr}$ be a Gaussian probability measure on $C^0(\overline{\Omega})$, and let $P : \theta \in C^0(\overline{\Omega}) \mapsto \sigma \in U$ be defined as $P_1$ or $P_2$ in Section 3.1. Then the function $\Phi^0 : C^0(\overline{\Omega}) \times \mathbb{R}^{JL} \to \mathbb{R}$ defined in (13), with $G^0 = F^0 \circ P : C^0(\overline{\Omega}) \to \mathbb{R}^{JL}$ satisfies assumptions 1-3 of Theorem 2.*

*Proof.* We have that $F^0 : U \to \mathbb{R}^{JL}$ is bounded, since from Corollary 1 we obtain

$$||F^0(\sigma)||_{C_e} \leq ||F^0(\sigma_U)||_{C_e} + ||F^0(\sigma) - F^0(\sigma_U)||_{C_e}$$

$$\leq \max\{1, ||F^0(\sigma_U)||_{C_e}\}(||\sigma||_{L^\infty(\Omega)} + c),$$

for some $\sigma_U \in U$ such that $||\sigma_U||_{L^\infty(\Omega)} = c$. In case $P = P_1$ then we have that

$$||G^0(\theta)||_{C_e} \leq C(||\exp(\theta)||_{L^\infty(\Omega)} + c) \tag{27}$$

$$\leq C(\exp(||\theta||_{L^\infty(\Omega)}) + c), \tag{28}$$

while if $P = P_2$, $||G^0(\theta)||_{C_e}$ is bounded by a constant since $P_2$ is uniformly bounded. Using the triangle inequality we have that

$$\Phi^0(\theta, z) \leq C(||z||^2_{C_e} + ||G^0(\theta)||^2_{C_e})$$

and therefore assumption 1 follows. To fulfill assumption 3 we note that we have

$$|\Phi^0(\theta, z_1) - \Phi^0(\theta, z_2)| = \frac{1}{2}|\langle z_1 + z_2 - 2G^0(\theta), z_1 - z_2\rangle_{C_e}|$$

$$\leq C(||z_1||_{C_e} + ||z_2||_{C_e} + 2||G^0(\theta)||_{C_e})||z_1 - z_2||_{C_e}$$

and the result follows. It remains to show that assumption 2 is also satisfied. Assume that $P$ is some map such that if $||\theta - \theta_n||_{L^\infty(\Omega)} \to 0$ then $P(\theta_n) \to P(\theta)$ uniformly or in measure. Then by Corollary 1 we have that $G^0 = F^0 \circ P$ is continuous at $\theta$. The continuity assumption on $P$ is true for $P = P_1$. By Proposition 1 it is true for $P = P_2$ at the points where the level sets have measure zero. However since we are assuming $\theta \sim \mu_{pr}$ and $\mu_{pr}$ is a Gaussian probability measure on $C^0(\overline{\Omega})$, it follows from Proposition 1 that $\theta$ has $\mu_{pr}$-almost surely this property and therefore assumption 3 is satisfied also in the case $P = P_2$. □

*Remark 3.* Using the same arguments, we can prove that when $P = P_1$ also for $G^\varepsilon(\theta)$ we have that

$$||G^\varepsilon(\theta)||_{C_e} \leq C(\exp(||\theta||_{L^\infty(\Omega)}) + c),$$

while when $P = P_2$, $||G^\varepsilon(\theta)||_{C_e}$ is bounded by a constant independent of $\theta$. Hence the posterior measure (8) based on the potential function $\Phi^\varepsilon$ is also well defined and Lipschitz contnuous in the data with respect to the Hellinger distance.

## 4.2 Modelling error and convergence analysis

Before moving to the numerical aspects of the problem, an investigation of the validity of our approach is necessary. First we observe that (6) can be rewritten as

$$z = F^0(\sigma^*) + \underbrace{(F^\varepsilon(\sigma^*) - F^0(\sigma^*))}_{m(\sigma^*)} + e, \quad e \sim \mathcal{N}(0, C_e). \tag{29}$$

The quantity $m(\sigma^*)$ represents the modelling error capturing the mismatch between the full multiscale model and the homogenized one. In particular (29) suggests that the observed data originating from the full multiscale model can be seen as data originating from the homogenized model, which are affected by two sources of errors: the noise and the modelling error. Both sources of errors can affect our predictions and we must take them into account when solving inverse problems to obtain good approximations of the unknown, especially when $\varepsilon$ is relatively large. For the modelling error we can show that we have in our case that $m(\sigma) \to 0$ as $\varepsilon \to 0$ for every $\sigma \in U$, as stated in the following theorem.

**Theorem 5.** *Let $\{A_\sigma^\varepsilon\}_{\varepsilon>0}$ be a sequence of matrices in $M(\alpha, \beta, \Omega)$ which G-converges to the matrix $A_\sigma^0 \in M(\alpha, \beta, \Omega)$, and let $m(\sigma) = \text{vec}(\{\tilde{m}_{jl}(\sigma)\}_{\substack{1 \leq j \leq J \\ 1 \leq l \leq L}})$, where*

$$\tilde{m}_{jl}(\sigma) = \int_{\Gamma_j} (\Lambda_{A_\sigma^\varepsilon} - \Lambda_{A_\sigma^0}) g_l \cdot \phi_j \, \mathrm{d}s, \quad j = 1, \ldots, J, \, l = 1, \ldots, L,$$

*where $\Gamma_j \subset \partial\Omega$ for all $j = 1, \ldots, L$, $\Gamma_i \cap \Gamma_j = \emptyset$ for $i \neq j$, and $\phi_j, g_l \in H^{1/2}(\partial\Omega)$ for all $j = 1, \ldots, J$, $supp(\phi_j) \subseteq \Gamma_j$ for all $j = 1, \ldots, J$. Then $||m(\sigma)||_{C_e} \to 0$ as $\varepsilon \to 0$.*

*Proof.* We have that for arbitrary $j$ and $l$ and $\forall \sigma \in U$

$$|\tilde{m}_{jl}(\sigma)| = \left| \int_{\Gamma_j} (\Lambda_{A_\sigma^\varepsilon} - \Lambda_{A_\sigma^0}) g_l \cdot \phi_j \, \mathrm{d}s \right|.$$

Since $supp(\phi_j) \subseteq \Gamma_j$, we have that

$$|\tilde{m}_{jl}(\sigma)| = \left| \int_{\partial\Omega} (\Lambda_{A_\sigma^\varepsilon} - \Lambda_{A_\sigma^0}) g_l \cdot \phi_j \, \mathrm{d}s \right|,$$

and using integration by parts

$$|\tilde{m}_{jl}(\sigma)| = \left| \int_\Omega (A_\sigma^\varepsilon \nabla u^\varepsilon - A_\sigma^0 \nabla u^0) \cdot \nabla \tilde{\phi}_j \, \mathrm{d}x \right|, \tag{30}$$

where $\tilde{\phi}_j$ is some function in $H^1(\Omega)$ whose trace is $\phi_j$. From G-convergence of $A_\sigma^\varepsilon$ to $A_\sigma^0$ we know that (30) converges to zero as $\varepsilon \to 0$. □

**Corollary 3.** Let $\mu^\varepsilon(\theta|z)$ and $\mu^0(\theta|z)$ be defined as in 8 and 14 respectively. Then we have that
$$d_{Hell}(\mu^0(\theta|z), \mu^\varepsilon(\theta|z)) \to 0 \text{ as } \varepsilon \to 0.$$

*Proof.* From the definition of the Hellinger distance we have that

$$2d_{Hell}^2(\mu^0, \mu^\varepsilon) = \int_{C^0(\overline{\Omega})} \left( \sqrt{\frac{d\mu^0}{d\mu_{pr}}} - \sqrt{\frac{d\mu^\varepsilon}{d\mu_{pr}}} \right)^2 \mu_{pr}(d\theta)$$

$$= \int_{C^0(\overline{\Omega})} \left( \frac{1}{\sqrt{C^0}} \exp\left(-\frac{1}{2}\Phi^0(\theta, z)\right) - \frac{1}{\sqrt{C^\varepsilon}} \exp\left(-\frac{1}{2}\Phi^\varepsilon(\theta, z)\right) \right)^2 \mu_{pr}(d\theta), \quad (31)$$

where $C^0$ and $C^\varepsilon$ are the two normalization constants such that $\mu^0(\theta|z)$ and $\mu^\varepsilon(\theta|z)$ are probability measures, i.e.,

$$C^0 = \int_{C^0(\overline{\Omega})} \exp(-\Phi^0(\theta, z))\mu_{pr}(d\theta), \quad C^\varepsilon = \int_{C^0(\overline{\Omega})} \exp(-\Phi^\varepsilon(\theta, z))\mu_{pr}(d\theta).$$

Let us notice that

$$|C^0 - C^\varepsilon| \leq \int_{C^0(\overline{\Omega})} \left| \exp(-\Phi^0(\theta, z)) - \exp(-\Phi^\varepsilon(\theta, z)) \right| \mu_{pr}(d\theta)$$

$$\leq \int_{C^0(\overline{\Omega})} \left| \Phi^0(\theta, z) - \Phi^\varepsilon(\theta, z) \right| \mu_{pr}(d\theta). \quad (32)$$

From (31) we get that

$$2d_{Hell}^2(\mu^0, \mu^\varepsilon) \leq I_1 + I_2,$$

where

$$I_1 = \frac{1}{C^0} \int_{C^0(\overline{\Omega})} \left( \exp\left(-\frac{1}{2}\Phi^0(\theta, z)\right) - \exp\left(-\frac{1}{2}\Phi^\varepsilon(\theta, z)\right) \right)^2 \mu_{pr}(d\theta),$$

$$I_2 = \left( \frac{1}{\sqrt{C^0}} - \frac{1}{\sqrt{C^\varepsilon}} \right)^2 C^\varepsilon.$$

We have that

$$I_1 \leq \frac{1}{4C^0} \int_{C^0(\overline{\Omega})} (\Phi^0(\theta, z) - \Phi^\varepsilon(\theta, z))^2 \mu_{pr}(d\theta),$$

and

$$I_2 \leq \frac{1}{4} \max\left\{ (C^0)^{-3}, (C^\varepsilon)^{-3} \right\} (C^0 - C^\varepsilon)^2$$

$$\leq C \int_{C^0(\overline{\Omega})} (\Phi^0(\theta, z) - \Phi^\varepsilon(\theta, z))^2 \mu_{pr}(d\theta),$$

where we have used (32). Using the definition of $\Phi^0$ and $\Phi^\varepsilon$ we find

$$2d_{Hell}^2(\mu^0, \mu^\varepsilon) \leq C \int_{C^0(\overline{\Omega})} (\Phi^0(\theta, z) - \Phi^\varepsilon(\theta, z))^2 \mu_{pr}(d\theta)$$

$$\leq C \int_{C^0(\overline{\Omega})} (2||z||_{C_e} + ||G^0(\theta)||_{C_e} + ||G^\varepsilon(\theta)||_{C_e})^2 ||G^0(\theta) - G^\varepsilon(\theta)||_{C_e}^2 \mu_{pr}(d\theta).$$

From Theorem 5 we have that $\lim_{\varepsilon \to 0} ||G^0(\theta) - G^\varepsilon(\theta)||_{C_e} = 0$. We also have that (see Corollary 2) if $P = P_1$ both $||G^0(\theta)||_{C_e}$ and $||G^\varepsilon(\theta)||_{C_e}$ are bounded by some scalar multiple of $(\exp(||\theta||_{L^\infty(\Omega)}) + 1)$, which is square integrable. Otherwise if $P = P_2$ both $||G^0(\theta)||_{C_e}$ and $||G^\varepsilon(\theta)||_{C_e}$ are bounded by a constant since $P_2$ is uniformly bounded, and again square integrability follows. Then by the Lebesgue's dominated convergence theorem it follows that $d_{Hell}(\mu^0, \mu^\varepsilon) \to 0$ as $\varepsilon \to 0$. □

The interpretation of the result is that when $\varepsilon$ is small we can neglect the modelling error, since it will be close to zero, and we do not need to take into account its probability distribution in the inversion process. However for larger values of $\varepsilon$, the mismatch between the observations and the data produced by the homogenized model might not be negligible, and using the coarse grained approach without taking into account the modelling error distribution may lead to bad predictions. In order to avoid that, we can correct the likelihood function, by approximating the probability distribution of the modelling error. We do so by using Algorithm 1 proposed in [9], which aims at approximating the mean $\overline{m}$ and the covariance $C_m$ of the modelling error distribution.

---
**Algorithm 1:** Approximation of modelling error distribution
**input** : prior measure, sample size $M$, map $P : \theta \mapsto \sigma$
**output:** mean $\overline{m}$ and covariance $C_m$ of the modelling error

1 Draw from the prior measure a sample of realizations $S = \{\theta_1, \ldots, \theta_M\}$
2 **for** $i = 1, \ldots, M$ **do**
3 $\quad m_i = G^\varepsilon(\theta_i) - G^0(\theta_i) = F^\varepsilon(P(\theta_i)) - F^0(P(\theta_i))$
4 **end**
5 $\overline{m} = \frac{1}{M} \sum_{i=1}^{M} m_i$
6 $C_m = \frac{1}{M} \sum_{i=1}^{M} (\overline{m} - m_i)(\overline{m} - m_i)^\top$

---

We assume a Gaussian distribution for the modelling error, so that $m \sim \mathcal{N}(\overline{m}, C_m)$ for all $\sigma$, and we can rewrite (29) as

$$z = F^0(\sigma^*) + e, \quad e \sim \mathcal{N}(\overline{m}, C_e + C_m). \tag{33}$$

Then as illustrated in Algorithm 1 the modelling error distribution is approximated offline. Only $M$ evaluations of the full multiscale model are needed. Hence, we use this approximation to modify the cost functional as in (34), and sample from the posterior by evaluating only the coarse homogenized model. We note that in (33) to apply the Bayesian framework for inverse problem, we still assume the independence of $e$ and $\theta$, despite the introduction of the modelling error in $e$. Nevertheless the practical usefulness of such algorithm has been shown in numerous works (see [6,9]). Then we may define the new likelihood as

$$\Phi^0(\theta, \overline{z}) = \frac{1}{2} ||\overline{z} - G^0(\theta)||^2_{C_e + C_m}, \tag{34}$$

where $\overline{z} = z - \overline{m}$. Note that conclusions about existence and well-posedness of the posterior measure are still valid under this definition of the potential function, which is equivalent to the one in (13), apart from the fact that observations $z$ are shifted by $\overline{m}$, and the covariance matrix is given by $C_e + C_m$.

## 5 Numerical approximation of the posterior density

The output of the Bayesian approach consists in the posterior measure. However in practice numerical sampling is needed to approximate the distribution, in order to obtain some meaningful informations (such as expected value or variance of the unknown, confidence intervals). As mentioned in Section 3.1 we consider as prior a Gaussian measure $\mu_{pr} = \mathcal{N}(\theta_{pr}, C_{pr})$ on $C^0(\overline{\Omega})$, and a map $P_i : C^0(\overline{\Omega}) \to U$, $i = 1, 2$, such that each draw from $\mu_{pr}$ is mapped into the admissible set $U$. In numerical experiments to reduce the dimensionality of the unknown we use a truncated Karhunen-Loève expansion. Indeed, each draw $\theta \sim \mathcal{N}(\theta_{pr}, C_{pr})$ can be represented as

$$\theta = \theta_{pr} + \sum_{k=1}^{\infty} \sqrt{\lambda_k} \xi_k \varphi_k \,,$$

where $\{\varphi_k, \lambda_k\}_{k=1}^{\infty}$ is an orthonormal set of eigenfunctions and eigenvalues of $C_{pr}$, and $\{\xi_k\}_{k=1}^{\infty}$ is an i.i.d. sequence with $\xi_1 \sim \mathcal{N}(0,1)$. We can then consider the truncated Karhunen-Loève expansion

$$\begin{aligned} \theta^K &= \theta_{pr} + \sum_{k=1}^{K} \sqrt{\lambda_k} \xi_k \varphi_k \\ &= \theta^K(\boldsymbol{\xi}) \,, \end{aligned} \tag{35}$$

where $\boldsymbol{\xi} = (\xi_1, \ldots, \xi_K)^\top$, and $\{\varphi_k, \lambda_k\}_{k=1}^{K}$ is the orthonormal set of eigenfunctions and eigenvalues of $C_{pr}$ corresponding to the $K$ largest eigenvalues. The unknown parameter is then parametrized by the $K$ coeffiecients $\{\xi_k\}_{k=1}^{K}$, which are a priori i.i.d. as $\mathcal{N}(0,1)$, and the inverse problem consists in approximating the posterior distribution of the $K$ coefficients by sampling from the posterior density $\pi^0(\boldsymbol{\xi}|z)$ which is given by

$$\begin{aligned} \pi^0(\boldsymbol{\xi}|z) &\propto \exp\left(-\frac{1}{2}\|z - G^0(\theta^K(\boldsymbol{\xi}))\|^2_{C_e} - \frac{1}{2}(\theta^K(\boldsymbol{\xi}) - \theta_{pr})^\top C_{pr}^{-1}(\theta^K(\boldsymbol{\xi}) - \theta_{pr})\right) \\ &= \exp\left(-\frac{1}{2}\|z - G^0(\theta^K(\boldsymbol{\xi}))\|^2_{C_e} - \frac{1}{2}\boldsymbol{\xi}^\top \boldsymbol{\xi}\right) \,. \end{aligned} \tag{36}$$

To sample from the posterior density we employ the Marcov chain Monte Carlo techniques (MCMC). Many algorithms belonging to the family of MCMC sampling methods are available in the literature. We decide to use the Metropolis Hastings (MH) method, which we illustrate just below. With this approach at each iteration we generate a new candidate $\boldsymbol{\eta} \in \mathbb{R}^K$ from a proposal density $q(\boldsymbol{\xi}^k, \boldsymbol{\eta})$, $q : \mathbb{R}^K \times \mathbb{R}^K \to \mathbb{R}^+$, where $\boldsymbol{\xi}^k$ is the current value of the variable. This new candidate is accepted with probability

$$a(\boldsymbol{\xi}^k, \boldsymbol{\eta}) = \min\left\{1, \frac{\pi^0(\boldsymbol{\eta}|z)q(\boldsymbol{\eta}, \boldsymbol{\xi}^k)}{\pi^0(\boldsymbol{\xi}^k|z)q(\boldsymbol{\xi}^k, \boldsymbol{\eta})}\right\} \,. \tag{37}$$

Otherwise the candidate is rejected and the chain remains at the current position $\boldsymbol{\xi}^k$. Note that if the proposal density is symmetric, i.e. $q(\boldsymbol{\xi}^k, \boldsymbol{\eta}) = q(\boldsymbol{\eta}, \boldsymbol{\xi}^k)$, (37) reduces to

$$a(\boldsymbol{\xi}^k, \boldsymbol{\eta}) = \min\left\{1, \frac{\pi^0(\boldsymbol{\eta}|z)}{\pi^0(\boldsymbol{\xi}^k|z)}\right\} \,.$$

In our experiments we consider the random walk proposal distribution to explore the density. Then

$$q(\boldsymbol{\xi}^k, \boldsymbol{\eta}) \propto \exp\left(-\frac{1}{2s^2}(\boldsymbol{\eta} - \boldsymbol{\xi}^k)^\top(\boldsymbol{\eta} - \boldsymbol{\xi}^k)\right) \,, \tag{38}$$

which is symmetric, and leads to Algorithm 2. The approximation of the target distribution improves as the numbers of samples $N_{sample}$ increases, and asymptotic convergence is guaranteed as $N_{sample} \to \infty$ under certain regularity properties of the target distribution and the proposal density. Therefore the results may be strongly dependent on the number of samples required, but also on the proposal density. In particular is a difficult task to establish when a sample is large enough. At the same time another general issue for MH is the choice of $s$ in (38), whose magnitude affects the speed at which the posterior distribution is explored and the number of rejected realizations.

## 6 Numerical approximation of the forward model

At each MH iteration we need to evaluate $\pi^0(\boldsymbol{\eta}|z)$ where $\boldsymbol{\eta}$ is a new candidate point, and $\pi^0$ is the posterior density given in (36). Hence, given $\sigma = P(\theta^K(\boldsymbol{\eta}))$, where $P : C^0(\overline{\Omega}) \to U$ is one of two maps introduced in Section 3.1, what remains to be discussed is how to evaluate $F^0(\sigma)$. The procedure can be substantially split in two parts:

1. Solve for each $l = 1, \ldots, L$
$$\begin{aligned} -\nabla \cdot (A^0_\sigma \nabla u^0) &= 0 &&\text{in } \Omega, \\ u^0 &= g_l &&\text{on } \partial\Omega, \end{aligned}$$
where $A^0_\sigma$ is the homogenized tensor corresponding to $A^\varepsilon_\sigma(x) = A(\sigma(x), x/\varepsilon)$.

2. Approximate the boundary fluxes, and evaluate (12).

**Algorithm 2:** Metropolis Hastings (MH)
___
**input** : posterior density $\pi^0(\boldsymbol{\xi}|z)$, starting point $\boldsymbol{\xi}^1 \in \mathbb{R}^K$, number of samples desired $N_{sample}$, symmetric proposal density $\mathcal{N}(0, s^2 I)$
**output:** sample of realization $S \in \mathbb{R}^{K \times N_{sample}}$
___
**1** Initialization: $k = 1$, $S = \boldsymbol{\xi}^1$
**2 while** $k < N_{sample}$ **do**
**3**   $\boldsymbol{\eta} = \boldsymbol{\xi}^k + sw$, $w \sim \mathcal{N}(0, I)$
**4**   $a(\boldsymbol{\xi}^k, \boldsymbol{\eta}) = \min\left\{1, \dfrac{\pi^0(\boldsymbol{\eta}|z)}{\pi^0(\boldsymbol{\xi}^k|z)}\right\}$
**5**   Draw $u \sim \mathcal{U}([0,1])$
**6**   **if** $a(\boldsymbol{\xi}^k, \boldsymbol{\eta}) > u$ **then**
**7**     accept: $\boldsymbol{\xi}^{k+1} = \boldsymbol{\eta}$
**8**   **else**
**9**     reject: $\boldsymbol{\xi}^{k+1} = \boldsymbol{\xi}^k$
**10**   **end**
**11**   $S = S \cup \boldsymbol{\xi}^{k+1}$, $k = k + 1$
**12 end**
___

### 6.1 Numerical homogenization

Since given $A^\varepsilon_\sigma$ the analytic solution for the corresponding $A^0_\sigma$ is often not known, numerical homogenization is needed. We employ the Finite Element Heterogeneous Multiscale Method (FE-HMM) which approximates the homogenized problem originating from (9) taking as input only the multiscale data. For a detailed analysis about FE-HMM we refer to [1, 2, 5]. We state here the simplest version involving only piecewise linear macro and micro simplicial elements. The method is based on a macro finite element space

$$S^1(\Omega, \mathcal{T}_H) = \{v^H \in H^1(\Omega) \,:\, v^H|_K \in \mathcal{P}^1(K), \,\forall K \in \mathcal{T}_H\},$$

where $\mathcal{T}_H$ is a partition of $\Omega$ in simplicial elements $K$ of diameter $H_K$, and $\mathcal{P}^1(K)$ is the space of linear polynomials on $K$. Then the homogenized tensor is approximated at each integration point by solving a micro problem on a sampling domain $K_\delta = x_K + (-\delta/2, \delta/2)^d$, $(\delta \geq \varepsilon)$. For a sampling domain $K_\delta$ we define a micro finite element space

$$S^1(K_\delta, \mathcal{T}_h) = \{v^h \in W(K_\delta) \,:\, v^h|_T \in \mathcal{P}^1(T) \,\forall T \in \mathcal{T}_h\}.$$

where

$$W(K_\delta) = W^1_{per} = \{v \in H^1_{per}(K_\delta); \int_{K_\delta} v \, dx = 0\}$$

in case we ask for periodic coupling, or

$$W(K_\delta) = H^1_0(K_\delta)$$

for a coupling with Dirichlet boundary conditions. Let $u^H$ be the approximate solution to problem (10). The numerical method reads as follows: find $u^H \in S^1(\Omega, \mathcal{T}_H)$, $u^H = g$ on $\partial \Omega$, such that
$$B_H(u^H, v^H) = F_H(v^H) \quad \forall v^H \in S_0^1(\Omega, \mathcal{T}_H),$$
where $S_0^1(\Omega, \mathcal{T}_H) = S^1(\Omega, \mathcal{T}_H) \cap H_0^1(\Omega)$, and
$$B_H(v^H, w^H) := \sum_{K \in \mathcal{T}_H} \frac{|K|}{|K_\delta|} \int_{K_\delta} A(\sigma(x_K), x/\varepsilon) \nabla v_K^h \cdot \nabla w_K^h \, dx, \tag{39}$$

and

$$F_H(v^H) := \sum_{K \in \mathcal{T}_H} |K|(fv^H)(x_K).$$

In (39) $v_K^h$ (respectively $w_K^h$) denotes the solution to the micro problem: find $v_K^h$ such that $v_K^h - v^H \in S^1(K_\delta, \mathcal{T}_h)$ and

$$\int_{K_\delta} A_\sigma^\varepsilon \nabla v_K^h \cdot \nabla z^h \, dx = 0 \quad \forall z^h \in S^1(K_\delta, \mathcal{T}_h). \tag{40}$$

## 6.2 Model order reduction

The main cost of the FE-HMM comes from the repeated solution of cell problems, whose number increases as we refine the macro mesh to obtain an appropriate approximation of the homogenized solution. This is particularly undesirable when solving inverse problems, since by using e.g. MH method one needs multiple evaluations of the cost functional for different realizations of the parameters of interest. We therefore combine the reduced basis methodology with the FE-HMM to drastically reduce the computational effort. For a detailed description and analysis of the method, called the Reduced Basis Finite Element Heterogeneous Multiscale Method (RB-FE-HMM), we mention [3]. The main idea is the following. During what is called the *offline stage* we select a small number $N$ of carefully precomputed micro solutions to construct a small subspace of microscopic functions which we call $S^N(Y)$. Then in the *online stage* we use these precomputed solutions to obtain fast evaluations of the homogenized tensor at the macroscopic quadrature points. The basis of $S^N(Y)$ are selected by using a greedy procedure. We randomly define a training set $\Xi_{Train} = \{(t_n, \eta_n) : 1 \le n \le N_{Train}, 1 \le \eta_n \le d\}$, where $t_n \in [\sigma^-, \sigma^+]$, while $\eta_n$ corresponds to the unit vector $\mathbf{e}_{\eta_n}$ of the canonical basis of $\mathbb{R}^d$. Let us note that if we map the domain $K_\delta$ into the reference domain $Y = (0, 1)^d$ through $x = G_{x_K}(y) = x_K + \delta(y - 1/2)$, we can define

$$B_H(v^H, w^H) := \sum_{K \in \mathcal{T}_H} |K| A^{0,h}(\sigma(x_K)) \nabla v^H(x_K) \cdot \nabla w^H(x_K), \tag{41}$$

with

$$(A_\sigma^{0,h}(x_K))_{ik} = \int_Y A_{x_K} (\nabla \chi_K^{i,h} + \mathbf{e_i}) \cdot (\nabla \chi_K^{k,h} + \mathbf{e_k}) \, dy, \tag{42}$$

where $A_{x_K} = A(\sigma(x_K), G_{x_K}(y))$, and $\chi_K^{i,h}$ (respectively $\chi_K^{k,h}$) is the solution of the micro problem

$$\begin{aligned} b(\chi_K^{i,h}, z^h) :&= \int_Y A_{x_K} \nabla \chi_K^{i,h} \cdot \nabla z^h \, dy \\ &= -\int_Y A_{x_K} \mathbf{e_i} \cdot \nabla z^h \, dy =: l_i(z^h) \; \forall z^h \in S^1(Y, \mathcal{T}_h). \end{aligned} \tag{43}$$

It can be shown [3] that

$$\frac{1}{|K_\delta|} \int_{K_\delta} A(\sigma(x_K), x/\varepsilon) \nabla v_K^h \cdot \nabla w_K^h \, dx = A^{0,h}(\sigma(x_K)) \nabla v^H(x_K) \cdot \nabla w^H(x_K) \,. \tag{44}$$

To select the micro problems in the offline stage, we start by computing the first basis function by solving the micro problem

$$\int_Y A(t_n, y) \nabla \zeta_1^h \cdot \nabla v^h \, dy = -\int_Y A(t_n, y) \mathbf{e}_{\eta_n} \cdot \nabla v^h \, dy = \forall v^h \in S^1(Y, \mathcal{T}_h) \,, \tag{45}$$

where $(t_n, \eta_n)$ is some point randomly drawn from $\Xi_{Train}$, and we initialize the space $S_1(Y) = \mathrm{span}\{\zeta_1^h\}$. We then successively we continue to add new basis functions to $S_1(Y)$ until convergence of an a posteriori error indicator below a certain tolerance. The output of the offline stage is then the reduced space

$$S_N(Y) = \mathrm{span}\{\zeta_1^h, \ldots, \zeta_N^h\} \,.$$

An efficient way to both implement the greedy algorithm and compute residuals for the a posteriori error control is crucial. A crucial assumption is that for a given $t_n \in [\sigma^-, \sigma^+]$, the tensor $A(t_n, y)$ is available in the affine form

$$A(t_n, y) = \sum_{q=1}^Q \Theta_q(t_n) A_q(y) \,, \forall y \in Y \,, \tag{46}$$

where $\Theta_q : [\sigma^-, \sigma^+] \to \mathbb{R}$, $q = 1, \ldots, Q$. In case the tensor is not directly available in the form (46), a greedy algorithm, called the empirical interpolation method (EIM), can be applied to obtain an affine approximation of the tensor [15]. Once $S_N(Y)$ has been constructed, we can then define a macro method similar to the FE-HMM. The method reads: find $u^{H,RB} \in S^1(\Omega, \mathcal{T}_H)$, $u^{H,RB} = g$ on $\partial \Omega$, such that

$$B_{H,RB}(u^{H,RB}, v^H) = \int_\Omega f v^H \, dx \quad \forall v^H \in S_0^1(\Omega, \mathcal{T}_H) \,, \tag{47}$$

where

$$B_{H,RB}(v^H, w^H) := \sum_{K \in \mathcal{T}_H} |K| A^{0,N}(\sigma(x_K)) \nabla v^H(x_K) \cdot w^H(x_K) \,, \tag{48}$$

and

$$(A^{0,N}(\sigma(x_K)))_{ik} = \int_Y A(\sigma(x_K), y)(\nabla \chi_K^{i,N} + \mathbf{e_i}) \cdot (\nabla \chi_K^{k,N} + \mathbf{e_k}) \, dy \,.$$

Here $\chi_K^{i,N}$ is the solution of (43) in the reduced basis space. Thanks to the affine representation of the tensor $A_\sigma^\varepsilon$, solving the micro problems in the reduced space consists with solving a small $N \times N$ linear system, which leads to a significant computational speed-up compared to classical numerical homogenization [3].

### 6.3 Approximate boundary flux computation

Once the approximated solution $u^{H,RB}$ has been computed, what is left is to evaluate boundary fluxes to obtain data defined through (12). To do so we employ a method based on a Galerkin projection, which we have introduced and analysed in [4]. For completeness, we briefly illustrate the algorithm. We assume the domain $\Omega$ to be a polygonal domain. We suppose that the approximate flux at the corners of $\Omega$ is specified from direct calculation with the given Dirichlet data (see [10] for details). We next introduce the following subspace of $S^1(\Omega, \mathcal{T}_H)$:

$$S_c^1(\Omega, \mathcal{T}_H) = \{v^H \in S^1(\Omega, \mathcal{T}_H) : v^H = 0 \text{ at the corners of } \Omega\} \,.$$

Let us denote as $S^1(\partial\Omega, \mathcal{T}_H)$, the finite dimensional space of functions which are restrictions on the boundary of functions in $S_c^1(\Omega, \mathcal{T}_H)$. Suppose we have computed the approximate solutions $u^{H,RB}$ and we want to approximate the boundary flux across $\partial\Omega$. The approximate flux can be constructed by assembling the function $\Lambda_{A_\sigma^{0,N}}^H g \in S^1(\partial\Omega, \mathcal{T}_H)$, such that

$$\int_{\partial\Omega} \Lambda_{A_\sigma^{0,N}}^H g \cdot v^H \, ds = B_{H,RB}(u^{H,RB}, v^H) - \int_\Omega f v^H \, dx \qquad (49)$$

$\forall v^H \in S_c^1(\Omega, \mathcal{T}_H)$, where

$$B_{H,RB}(v^H, w^H) = \sum_{K \in \mathcal{T}_H} |K| A^{0,N}(\sigma(x_K)) \nabla v^H \cdot \nabla w^H \,.$$

Let us remark that $u^{H,RB}$ has been already computed, and so constructing $\Lambda_{A_\sigma^{0,N}}^H g$ leads to solving a linear systems whose unknowns are the values of the boundary flux at the nodes of $\partial\Omega$.

## 6.4 Summary of the numerical procedure to solve the multiscale inverse problem

The numerical scheme for solving the multiscale inverse problem given the perturbed observations $z \in \mathbb{R}^{JL}$ can be then summarized as follows.

1. Compute in an offline stage a reduced space of precomputed microscopic functions as described in Section 6.

2. Compute in an offline stage the set $\{\varphi_k, \lambda_k\}_{k=1}^K$ of eigenfunctions and eigenvalues of the prior covariance $C_{pr}$, so that for a point $\boldsymbol{\xi} \in \mathbb{R}^K$ we have that

$$\theta^K(\boldsymbol{\xi}) = \theta_{pr} + \sum_{k=1}^K \sqrt{\lambda_k} \xi_k \varphi_k \,,$$

   and $\sigma^K(\boldsymbol{\xi}) = P(\theta^K(\boldsymbol{\xi}))$, where $P = P_1$ or $P = P_2$.

3. Sample online from the posterior distribution. In particular for a new realization $\boldsymbol{\eta}$, in order to evaluate $\pi^0(\boldsymbol{\eta}|z)$, do the following for each $l = 1, \ldots, L$.

   - Solve (47), (48) for $\sigma = \sigma^K(\boldsymbol{\eta})$.
   - Approximate the corresponding flux at the boundary by solving (49) and evaluate (12).

# 7 Numerical experiments

In this section we will present some numerical experiments to illustrate our multiscale Bayesian algorithm for inverse problems. We start by explaining how observed data are collected. We then solve the inverse problem for different macroscopic parametrizations. At first, we consider an affine parametrization of the form $A^\varepsilon(x) = \sigma^*(x)B^\varepsilon(x) = \sigma^*(x)B(x/\varepsilon)$, so that the function $\sigma^*$ controls the amplitude of the characteristic micro oscillations. Let us point out that for this choice we have that $A^0(x) = \sigma^*(x)B^0$, and thus the use of reduced basis methods for solving the forward problem is not required. This simple problem allows us to perform numerous tests to quantify the sensitivity of the method with respect to the several parameters involved in the approximation, such as $\varepsilon$, the size of the microscopic oscillations, $K$, the number of terms in the truncated Karhunen-Loève expansion, and $L$, the number Dirichlet data. Then we will consider two different non-affine macroscopic parametrizations, one controlling the orientation of the micro oscillations, the other the volume fraction of a fictious layered material. For these problems we make the following choice of parameters for the RB-FE-HMM offline stage: $h/\varepsilon = 1/64$, $\delta = \varepsilon$, $tol_{RB} = 10^{-11}$, where $tol_{RB}$ is a prescribed tolerance used as stopping criterion for the greedy algorithm employed to select the reduced basis functions.

## 7.1 Set-up

The computational domain is the unit square

$$\Omega = \{x = (x_1, x_2) \,:\, 0 < x_1, x_2 < 1\}\,.$$

We approximate the solution to problem (1) by means of the Finite Element Method (FEMs) using a very fine discretization $h_{obs} << \varepsilon$. The forward homogenized problem is instead computed using a macro mesh size $H = 1/64$. The problem is solved for different Dirchlet conditions $\{g_l\}_{l=1}^L$. In particular we take $\{g_l\}_{l=1}^L = \{\sqrt{\lambda_l}\varphi_l\}_{l=1}^L$, where $\{(\lambda_l, \varphi_l)\}_{l=1}^L$ are the $L$ eigenpairs corresponding to the $L$ smallest eigenvalues associated to the one dimensional discrete Laplacian operator. Each $g_l$ is then projected on the boundary $\partial\Omega$ to define the corresponding Dirichlet condition. This procedure ensures that the functions $\{g_l\}_{l=1}^L$ are smooth and orthonormal, so that each experiment contributes differently one from another. Moreover $\|\nabla g_l\|_{L^2(\partial\Omega)} < C$, where $C$ is a constant independent of $L$. Finally, we consider $J = 12$ boundary portions $\Gamma_j \subset \partial\Omega$, three for each side of the computational domain as shown in Figure 1. Each $\Gamma_j$ has length equal to 0.2. The functions $\phi_j$ appearing in (5) are hat functions with $supp(\phi_j) = \Gamma_j$ which take value one at the midpoint of each $\Gamma_j$. Once the observed data have been computed, they are perturbated by the noise given by $e = 10^{-4}w$, $w \sim \mathcal{N}(0, I)$. Let $p_i$ and $p_j$ two nodes of the macro triangulation $\mathcal{T}_H$, and let $N_H$ the total number of nodes defining $\mathcal{T}_H$. Note that $N_H = H^{-2}$. The covariance matrix in the prior measure $\mu_{pr} = \mathcal{N}(\theta_{pr}, C_{pr})$ is then $C_{pr} \in \mathbb{R}^{N_H \times N_H}$ defined as

$$(C_{pr})_{ij} = \gamma \exp\left(-\frac{\|p_i - p_j\|}{\lambda}\right), \gamma, \lambda \in \mathbb{R}^+, \tag{50}$$

while the prior mean $\theta_{pr}$ is some function in $C^0(\overline{\Omega})$. We set different values for $\gamma, \lambda$ and $\theta_{pr}$ depending on the macroscopic parametrization we want to retrieve. In particular $\lambda > 0$ is a correlation length that describes how the values at different positions of the functions supported by the prior measure are related, while $\gamma > 0$ is the amplitude scaling factor.

## 7.2 2D affine parametrization (amplitude of oscillations)

In this first set of numerical experiments we consider the tensor $A^\varepsilon_{\sigma^*}$ given by

$$a_{11}(\sigma^*(x), x/\varepsilon) = \sigma^*(x)\left(\cos^2\left(\frac{2\pi x_1}{\varepsilon}\right) + 1\right),$$
$$a_{22}(\sigma^*(x), x/\varepsilon) = \sigma^*(x)\left(\sin\left(\frac{2\pi x_2}{\varepsilon}\right) + 2\right),$$
$$a_{12}(\sigma^*(x), x/\varepsilon) = a_{21}(\sigma^*(x), x/\varepsilon) = 0\,,$$

where

$$\sigma^*(x) = 1.3 + 0.3\mathbb{1}_{D_1} - 0.4\mathbb{1}_{D_2}\,,$$

and

$$D_1 = \{x = (x_1, x_2) \,:\, (x_1 - 5/16)^2 + (x_2 - 11/16)^2 \leq 0.025\}\,,$$
$$D_2 = \{x = (x_1, x_2) \,:\, (x_1 - 11/16)^2 + (x_2 - 5/16)^2 \leq 0.025\}\,.$$

The task of the problem is to retrieve the function $\sigma^*$, which is shown together with the component $a_{11}^\varepsilon$ of the tensor, $\varepsilon = 1/64$, in Figure 2.

**Sensitivity with respect to $\varepsilon$**

We start by studying how different choices of $\varepsilon$ can affect our predictions. The computations are reported in Figure 4. We briefly describe the setting. We compute numerically by means of a resolved FEM synthetic observations for different values of $\varepsilon = \{1/4, 1/8, 1/16, 1/8, 1/64\}$, for $L = 6$

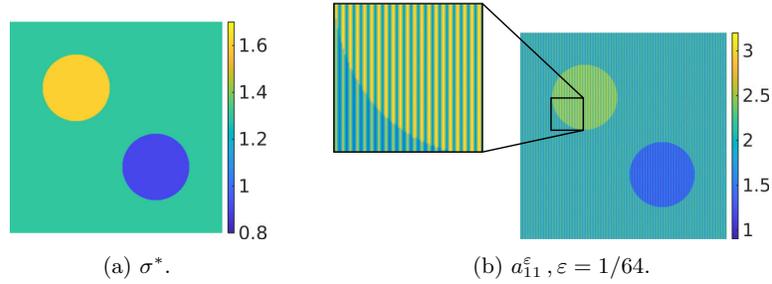

(a) $\sigma^*$.     (b) $a_{11}^\varepsilon$, $\varepsilon = 1/64$.

Figure 2: Representation of the true spatial field $\sigma^*$ and the first component of the highly oscillating tensor for Problem 7.2.

different Dirichlet conditions. We consider a truncated Karhunen-Loève expansion with $K = 60$. The prior measure $\mu_{pr}$ on $\theta \in C^0(\overline{\Omega})$ is $\mathcal{N}(\theta_{pr}, C_{pr})$, with $\theta_{pr} = \log 1.3$ and $C_{pr}$ defined in (50) with $\gamma = 0.05$ and $\lambda = 0.5$. In particular the choice of $\theta_{pr} = \log 1.3$ is such that the resulting log-normal distribution on the admissible set $U$ has median 1.3. For $\mu_{pr}$ measure on $C^0(\overline{\Omega})$ and $P : C^0(\overline{\Omega}) \to U$, we denote as $P^\# \mu_{pr}$ the pushed forward prior on the admissible set $U$ under $P$. We then push each draw $\theta$ into the admissible set through the function $P_1 : \theta \mapsto \exp(\theta)$. Example of realizations from the pushed forward prior $P_1^\# \mu_{pr}$ are show in Figure 3. We draw then $2 \times 10^5$ samples from the posterior distribution (36) using MH. The parameters $s$ is set to 0.01. The starting point is $\boldsymbol{\xi}^1 = 0 \in \mathbb{R}^K$. With this choice of the parameters we obtain an acceptance rate of about 27% for all choices of $\varepsilon$. In Figure 4 we plot for each $\varepsilon$ the quantities $P_1(\mathbb{E}[\theta^K(\boldsymbol{\xi})])$, $\mathbb{E}[P_1(\theta^K(\boldsymbol{\xi}))]$, and the variance $\mathrm{Var}[P_1(\theta^K(\boldsymbol{\xi}))]$. The first quantity is produced by computing first the mean on the Banach space $C^0(\overline{\Omega})$ and then pushing it into the admissible set $U$ through $P_1 : C^0(\overline{\Omega}) \to U$. Moreover we also show the approximation of the posterior density for the first three coefficients in the truncated Karhunen-Loève expansion. We notice that with $\varepsilon = 1/4$ we get inaccurate predictions about the quantity of interest, while already with $\varepsilon = 1/8$ the approximation of the posterior mean is in good agreement with Figure 2. The source of error for large $\varepsilon$ comes from the discrepancy between the multiscale model from where the observations are obtained and the homogenized model used for solving the inverse problem.

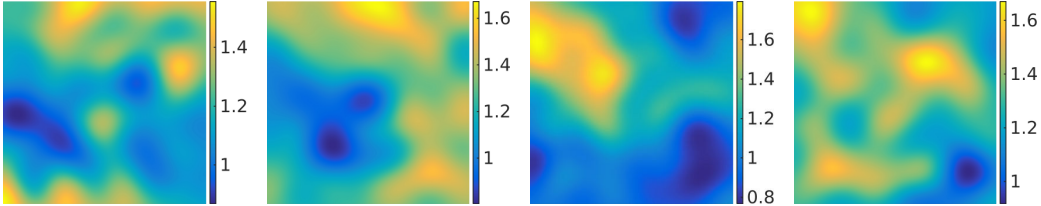

Figure 3: Four samples from the prior density used in Problem 7.2.

**Approximation of the modelling error distribution**

As seen in Figure 4 for large values of $\varepsilon$ the modelling error (the discrepancy between the fine scale and the homogenized problems) pollutes the posterior prediction. Therefore, we perform again the same experiment for $\varepsilon = 1/4$, but taking into account the modelling error as described in Section 4.2. We approximate the modelling error distribution, by computing its mean and covariance using Algorithm 1 and include these quantities into the posterior density definition according to 34. We perform the experiment for various number of sample sizes $M$ to approximate the modelling error distribution, namely $M = \{5, 10, 20\}$. The parameters such as $K$ and $L$ are identical to the previous numerical test. Numerical results are shown in Figure 5. In particular we can observe how already with $M = 5$ we can manage to significantly improve the results reported in Figure 4 for $\varepsilon = 1/4$.

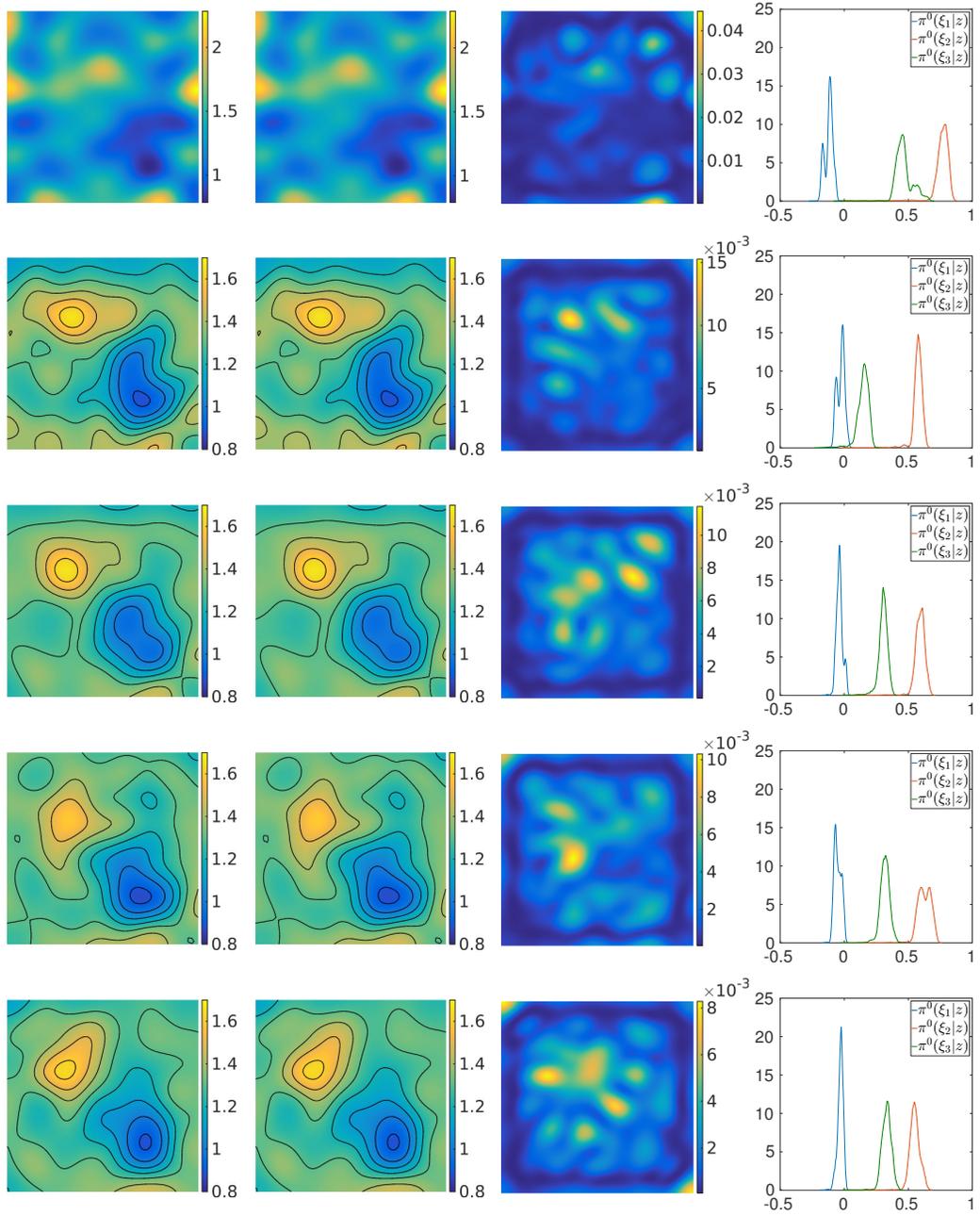

Figure 4: Comparison of numerical approximations of the posterior density for Problem 7.2, obtained with different values of $\varepsilon$. From left to right the plotted quantities are $P_1(\mathbb{E}[\theta])$, $\mathbb{E}[P_1(\theta)]$, $\text{Var}[P_1(\theta)]$, and the posterior density of the three first coefficients of the truncated Karhunen-Loève expansion, corresponding to $\varepsilon = \{1/4, 1/8, 1/16, 1/32, 1/64\}$. The length scale $\varepsilon$ decreases from the top to the bottom. The other parameters are $H = 1/64$, $L = 6$, $K = 60$.

**Sensitivity with respect to $L$ (number of Dirichlet data)**

Next we investigate the sensitivity of the approximated solution with respect to the parameter $L$, denoting the number of different Dirichlet conditions used to produce the observations. The setting is the same as in the previous numerical experiments, except that $\varepsilon$ is fixed and equal to $1/64$, while $L = \{2, 4, 6\}$. Numerical results are shown in Figure 6. We notice that for $L = 2$ the variance is significantly larger than for $L = 4$ or $L = 6$, which indicates more uncertainty about the approximated solution. This is also visible from the approximation of the posterior density obtained for the three first coefficients of the Karhunen-Loève expansion.

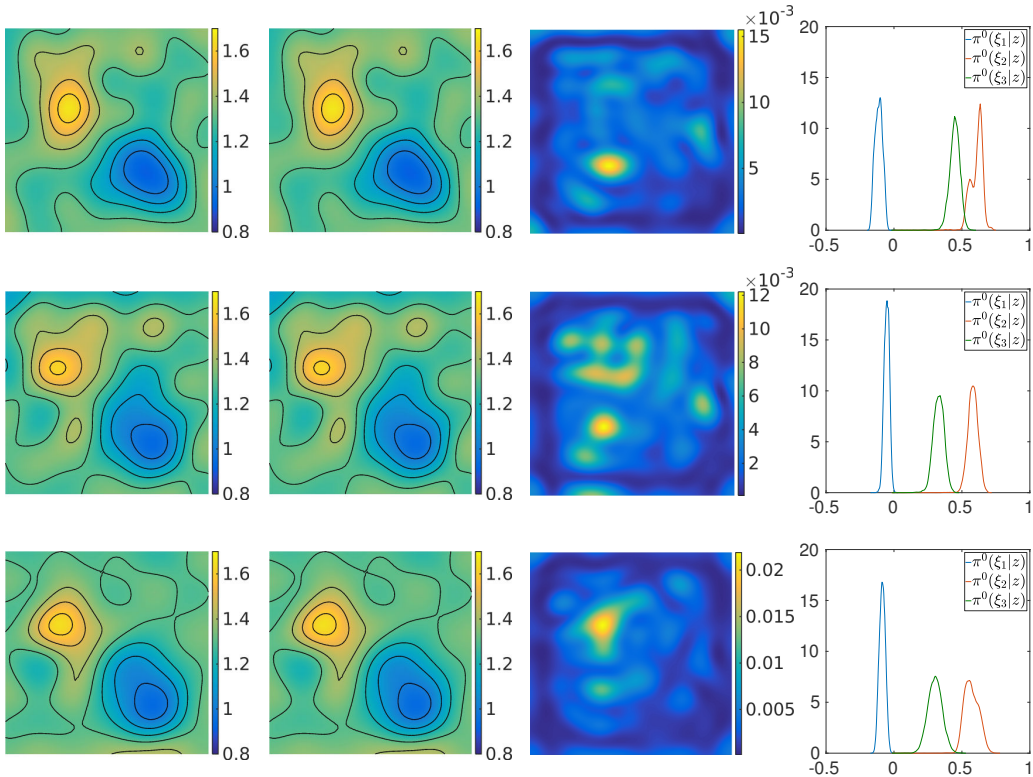

Figure 5: Comparison of numerical approximations of the posterior density for Problem 7.2 obtained with $\varepsilon = 1/4$, for different values of $M$, the sample size used to approximate the modelling error distribution. From left to right the plotted quantities are $P_1(\mathbb{E}[\theta])$, $\mathbb{E}[P_1(\theta)]$, $\mathrm{Var}[P_1(\theta)]$, and the posterior density of the three first coefficients of the truncated Karhunen-Loève expansion. The value of $M$ is 5 in the first row, 10 in the second one, and 20 in the third row. The other parameters are $H = 1/64$, $L = 6$, $K = 60$.

**Sensitivity with respect to $K$ (number of terms in the truncated KL expansion)**

Finally we examine how the size of the truncated Karhunen-Loève expansion affects our predictions. We perform experiments for $K = \{10, 20, 30, 40, 50, 60\}$, while $L$ and $\varepsilon$ are fixed, set to 6 and $1/64$ respectively. In particular we mention that for smaller $K$, a coarser mesh can be used for the forward discrete problem, leading to a significant saving of the computational cost. The results are shown in Figure 7. We can observe how the lowest Karhunen-Loève modes are able to determine the main geometric structure of the parameter of interest, while by increasing the number of eigenvalues/eigenfunctions we obtain a better sampling of the quantity of interest. This can be noticed from the plot of the posterior mean and variance we obtain with $K = \{40, 50, 60\}$. Such a result suggests the possibility of investigating the implementation of a Metropolis-Hastings algorithm on multiple levels, with an approximation of the distribution of the lowest modes on a coarse mesh, while performing fewer samples for the highest modes on a finer mesh to guarantee a proper sample of the posterior density.

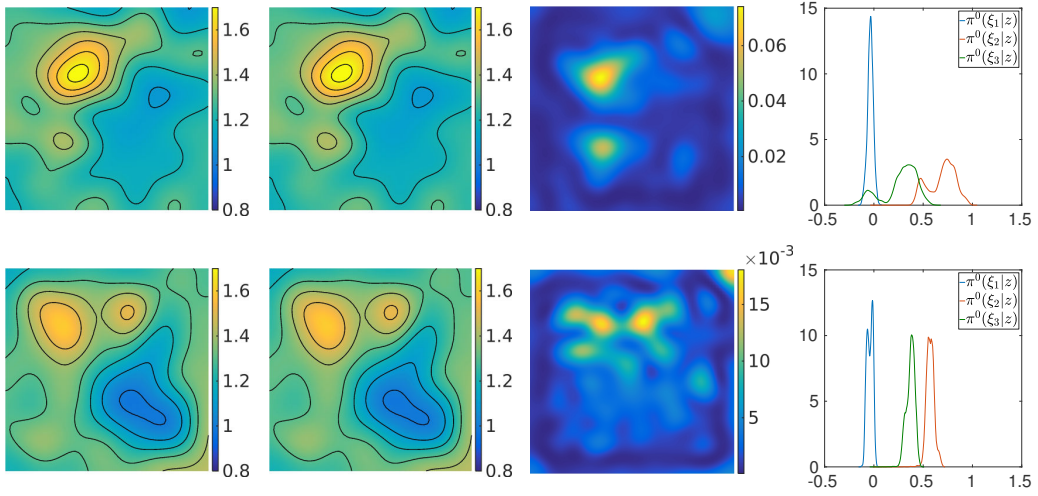

Figure 6: Comparison of numerical approximations of the posterior density for Problem 7.2 obtained for different values of $L$, the number of Dirchlet data. From left to right the plotted quantities are $P_1(\mathbb{E}[\theta])$, $\mathbb{E}[P_1(\theta)]$, $\text{Var}[P_1(\theta)]$, and the posterior density of the three first coefficients of the truncated Karhunen-Loève expansion. In the first row $L = 2$, in the second one $L = 4$. For $L = 6$ see last row in Figure 4. The other parameters are $H = 1/64$, $\varepsilon = 1/64$, $K = 60$.

### 7.3 2D non-affine parametrization (orientation of oscillations)

Now we consider the case where the function $\sigma^*$ controls the angle of the oscillations which characterize the full tensor $A^\varepsilon_{\sigma^*}$. The tensor is defined as

$$\begin{aligned}
a_{11}(\sigma^*(x), x/\varepsilon) &= \sin\left(\frac{4\pi \mathbf{e_1}^\top Q x}{\varepsilon}\right) + 1.5\,, \\
a_{22}(\sigma^*(x), x/\varepsilon) &= \cos^2\left(\frac{2\pi \mathbf{e_2}^\top Q x}{\varepsilon}\right) + 1\,, \\
a_{12}(\sigma^*(x), x/\varepsilon) &= a_{21}(\sigma^*(x), x/\varepsilon) = 0\,,
\end{aligned} \qquad (51)$$

where $Q = Q(\sigma^*(x))$ is a rotation matrix depending on $\sigma^* : \Omega \to \mathbb{R}$

$$Q(\sigma^*(x)) = \begin{pmatrix} \cos(2\pi\sigma^*(x)) & \sin(2\pi\sigma^*(x)) \\ -\sin(2\pi\sigma^*(x)) & \cos(2\pi\sigma^*(x)) \end{pmatrix}, \qquad (52)$$

and

$$\sigma^*(x) = a + b \mathbb{1}_D\,, \quad D \subset \Omega\,, a, b \in \mathbb{R}\,.$$

We consider the case where $D$ is the circle defined as

$$D = \{x = (x_1, x_2) : (x_1 - 1/3)^2 + (x_2 - 1/3)^2 \leq 0.05\}\,.$$

In Figure 8 we show the function $\sigma^*$ and the first component of the tensor $a^\varepsilon_{11}$. From (51)-(52) it can be observed that different values of $a$, $b$ for $\sigma^*$ can lead to the same rotation of the oscillations, and in general to the same tensor $A^\varepsilon_{\sigma^*}$. To ensure uniqueness we assume to know a priori the values of $a$, $b$. We take $a = 1$ and $b = 0.25$. Our task is thus to recover the region $D \subset \Omega$. To do so we consider as admissible set for the unknown the one defined through the level set function $P_2 : C^0(\overline{\Omega}) \to U$ introduced in Section 3.1. The prior measure on $C^0(\overline{\Omega})$ is defined as in 50 with $\theta_{pr} = 1$, $\gamma = 0.025$, and $\lambda = 0.5$. Each draw from $\mu_{pr}$ is then mapped into $U$ through $P_2$ defined as

$$P_2(\theta) = 1 \mathbb{1}_{\Omega_1} + 1.25 \mathbb{1}_{\Omega_2}\,,$$

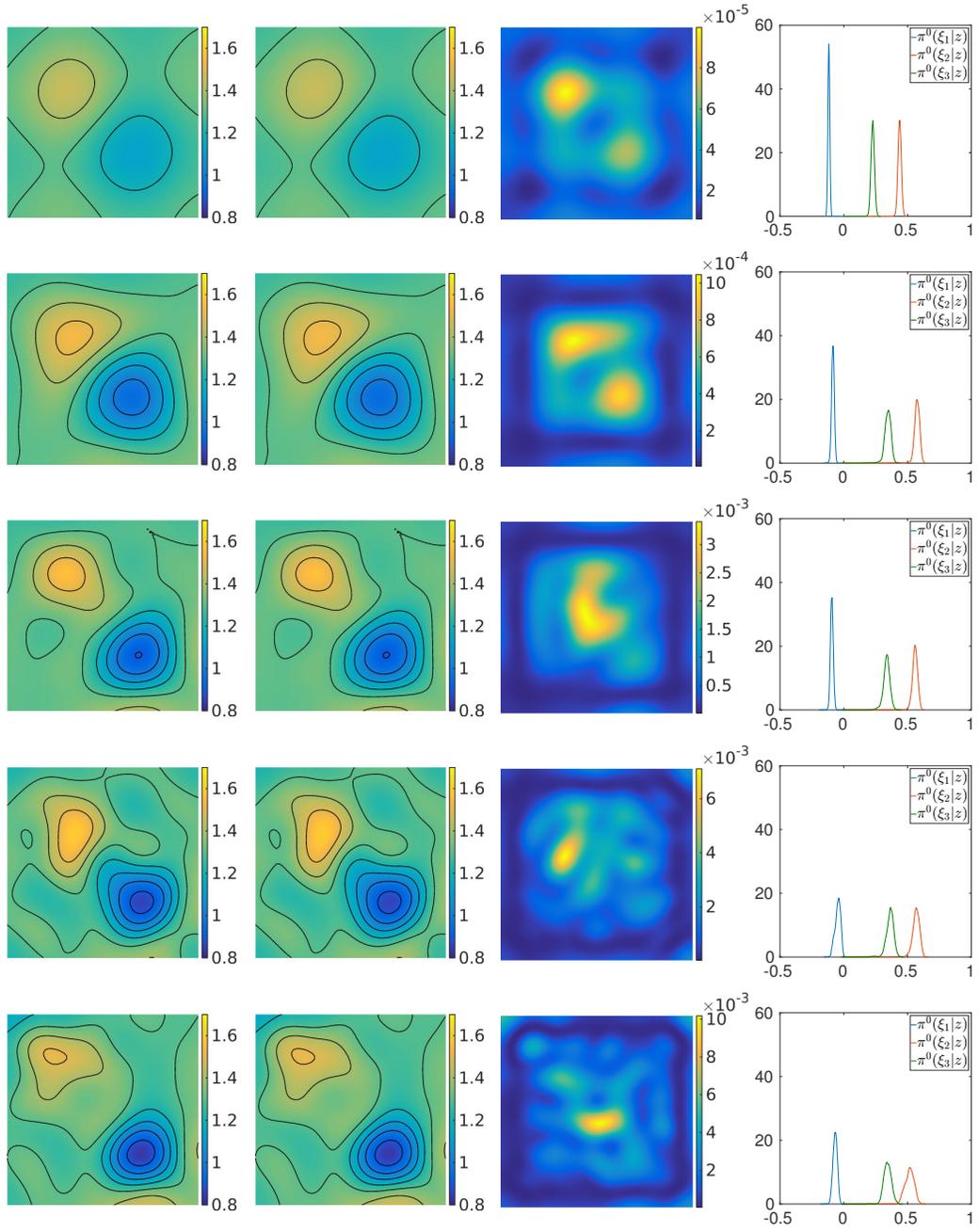

Figure 7: Comparison of numerical approximations of the posterior density for Problem 7.2 obtained for different values of $K$, the number of coefficients in the truncated Karhunen-Loève expansion. From left to right the plotted quantities are $P_1(\mathbb{E}[\theta])$, $\mathbb{E}[P_1(\theta)]$, $\text{Var}[P_1(\theta)]$, and the posterior density of the three first coefficients of the truncated Karhunen-Loève expansion, corresponding to $K = \{10, 20, 30, 40, 50\}$. The parameter $K$ increases from the top to the bottom. For $K = 60$ see last row in Figure 4. The other parameters are $H = 1/64$, $\varepsilon = 1/64$, $L = 6$.

where
$$\Omega_1 = \{x \in \Omega : -\infty < \theta(x) \leq 1\}, \quad \Omega_2 = \{x \in \Omega : 1 < \theta(x) < \infty\},$$

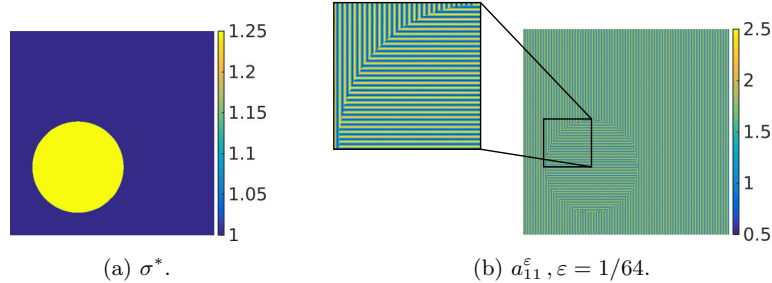

(a) $\sigma^*$.

(b) $a_{11}^\varepsilon$, $\varepsilon = 1/64$.

Figure 8: Representation of the true spatial field $\sigma^*$ and the first component of the highly oscillating tensor for the non-affine case Problem 7.3 (orientation of oscillations).

so that $\overline{\Omega} = \Omega_1 \cup \Omega_2$, $\Omega_1 \cap \Omega_2 = \emptyset$. Four examples of draws from the pushed forward prior $P_2^\# \mu_{pr}$ are reported in Figure 9. We obtain data for $\varepsilon = 1/64$ and approximate the modelling error distribution by using Algorithm 1 with $M = 20$. The parameters $K$ and $L$ are set to 60 and 6 respectively. Then we approximate the posterior by using Metropolis-Hastings by drawing $4 \times 10^5$ samples using $s = 0.02$. For this choice of the parameters we get an acceptance ratio during the sampling of about 73%. In Figure 10 we plot the quantities $P_2(\mathbb{E}[\theta])$, $\mathbb{E}[P_2(\theta)]$, and $\text{Var}[P_2(\theta)]$. In particular $P_2(\mathbb{E}[\theta])$ preserves the binary field property of the admissible set, while the estimate $\mathbb{E}[P_2(\theta)]$ gives a better understanding of the uncertainty across the interface where the discontinuity takes place. This uncertainty is also reflected by the plot of the variance $\text{Var}[P_2(\theta)]$. The numerical results show good agreement with Figure 8.

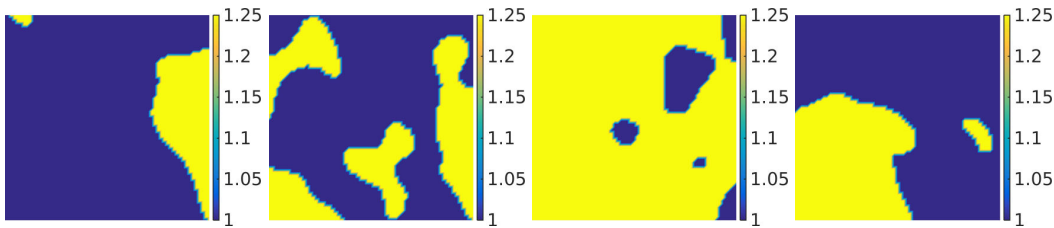

Figure 9: Four samples from the level set prior used in Problem 7.3.

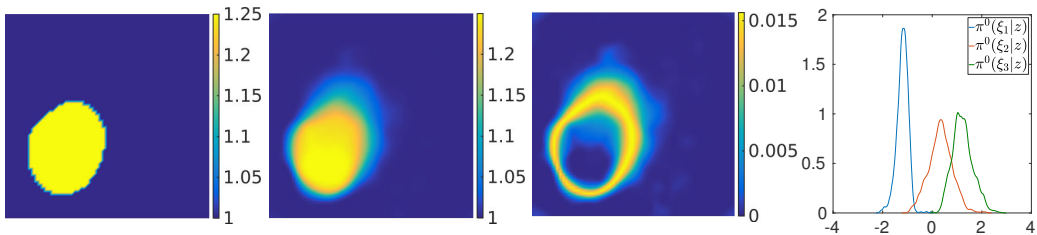

Figure 10: Numerical results for the non-affine parametrization Problem 7.3 (orientation of oscillations). From left to right the plotted quantities are $P_2(\mathbb{E}[\theta])$, $\mathbb{E}[P_2(\theta)]$, $\text{Var}[P_2(\theta)]$, and the posterior density of the three first coefficients of the truncated Karhunen-Loève expansion. The values of the parameters are $H = 1/64$, $\varepsilon = 1/64$, $M = 20$, $L = 6$, $K = 60$.

## 7.4 2D non-affine parametrization (volume fraction)

We conclude the numerical experiments by considering the case where $A^\varepsilon$ represents the conductivity of a hypothetical two phase layered material. In this case the macroscopic function $\sigma^* : \Omega \to [0,1]$ determines the volume fraction of each component. Then tensor is defined as

$$a_{11}(\sigma^*(x), x/\varepsilon) = a_{22}(\sigma^*(x), x/\varepsilon) = \begin{cases} 2 \text{ if } 0 \leq (x_2 \text{ mod } \varepsilon)/\varepsilon < \sigma^*(x) \\ 1 \text{ if } \sigma^*(x) \leq (x_2 \text{ mod } \varepsilon)/\varepsilon < 1 \end{cases}, \quad (53)$$

$$a_{12}(\sigma^*(x), x/\varepsilon) = a_{21}(\sigma^*(x), x/\varepsilon) = 0.$$

We consider the case where $\sigma^*$ is defined as

$$\sigma^*(x) = \sum_{i=1}^{n} c_i \mathbb{1}_{D_i}, \quad D_i \subset \Omega, c_i \in [0,1],$$

$D_i \cap D_j = \emptyset$ for $i \neq j$, $\cup_{i=1}^n D_i = \Omega$. Again we assume to know a priori the values $\{c_i\}_{i=1}^n$ that the function $\sigma^*$ can take, and our goal is to recover the different regions $\{D_i\}_{i=1}^n$. We note that knowing the range of possible values for $\sigma^*$ allows us to efficiently use the RB method and in particulaer the EIM algorithm. For our problem we set $n = 4$, $c_1 = 0.8, c_2 = 0.6, c_3 = 0.4, c_4 = 0.2$, and we make the following choice for the sets $\{D_i\}_{i=1}^4$

$$D_1 = \{x = (x_1, x_2) : 0 \leq x_1 \leq 0.25\},$$
$$D_2 = \{x = (x_1, x_2) : 0.25 < x_1 \leq 0.5\},$$
$$D_3 = \{x = (x_1, x_2) : 0.5 < x_1 \leq 0.75\},$$
$$D_4 = \{x = (x_1, x_2) : 0.75 < x_1 \leq 1\}.$$

The true field $\sigma^*$ and the first component of the multiscale tensor are shown in Figure 11. We

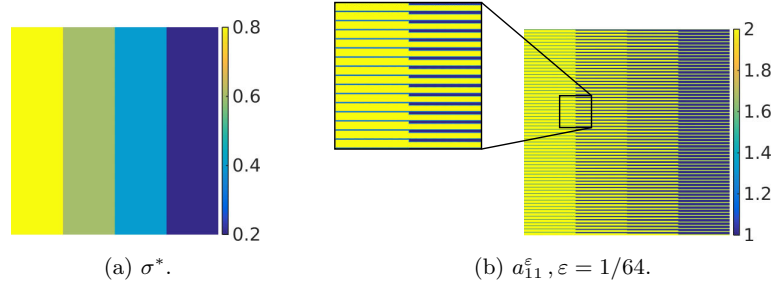

(a) $\sigma^*$.      (b) $a_{11}^\varepsilon$, $\varepsilon = 1/64$.

Figure 11: Representation of the true spatial field $\sigma^*$ and the first component of the highly oscillating tensor for the non-affine case Problem 7.4 (volume fraction).

consider for this last numerical experiment a macro discretization with mesh size $H = 1/32$, and a Gaussian prior measure $\mu_{pr}$ on $C^0(\overline{\Omega})$ as for the previous numerical tests, with $\theta_{pr} = 0.5$, $\gamma = 0.05$, $\lambda = 0.5$. The function $P_2 : C^0(\overline{\Omega}) \to U$, is instead defined as

$$P_2(\theta) = c_1 \mathbb{1}_{\Omega_1} + c_2 \mathbb{1}_{\Omega_2} + c_3 \mathbb{1}_{\Omega_3} + c_4 \mathbb{1}_{\Omega_4},$$

where

$$\Omega_1 = \{x \in \Omega : 0.6 < \theta(x) < \infty\},$$
$$\Omega_2 = \{x \in \Omega : 0.4 < \theta(x) \leq 0.6\},$$
$$\Omega_3 = \{x \in \Omega : 0.2 < \theta(x) \leq 0.4\},$$
$$\Omega_4 = \{x \in \Omega : -\infty < \theta(x) \leq 0.2\}.$$

Four samples from the prior $P_2^\# \mu_{pr}$ are shown in Figure 12. To solve the problem the observations are obtained for $\varepsilon = 1/64$. The modelling error distribution is approximated offline by using Algorithm 1 with $M = 20$. The parameter $K$ and $L$ are set to 60 and 6 respectively. We draw $4 \times 10^5$ samples from the posterior distribution using Algorithm 2 setting $s = 0.01$, leading to an acceptance ratio of 44%. The numerical results are shown in Figure 13, and are in good agreement with Figure 11.

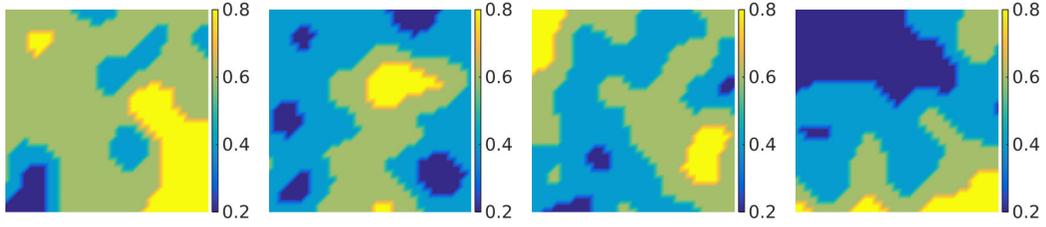

Figure 12: Four samples from the level set prior used in Problem 7.4.

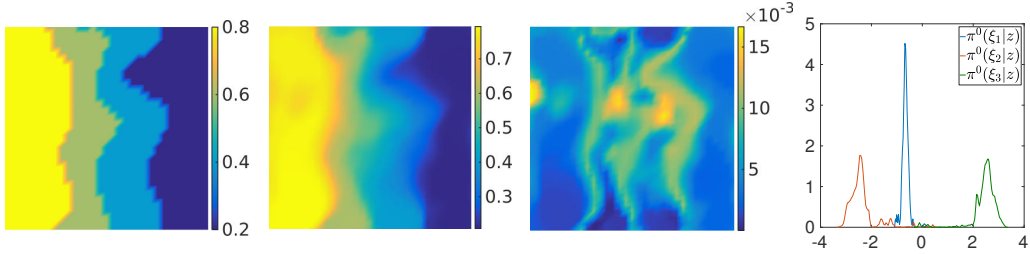

Figure 13: Numerical results for the non-affine parametrization Problem 7.4 (volume fraction). From left to right the plotted quantities are $P_2(\mathbb{E}[\theta])$, $\mathbb{E}[P_2(\theta)]$, $\mathrm{Var}[P_2(\theta)]$, and the posterior density of the three first coefficients of the truncated Karhunen-Loève expansion. $H = 1/32$, $\varepsilon = 1/64$, $M = 20$, $L = 6$, $K = 60$

## Conclusion

We have presented a new strategy for solving Bayesian multiscale inverse problems based on numerical homogenization and model order reduction. Our method allows to recover the full fine scale tensor under the assumption that the microscopic structure of the fine scale tensor is known to us but its macroscopic behaviour is unknown. Practical examples include multi-phase mediums, whose constituents are known, but their respective volume fraction or macroscopic orientation are unknown. We then proved the existence and well-posedness of the effective posterior measure obtained by homogenization of the forward operator. By means of G-convergence we showed that the fine scale posterior measure converges to the homogenized posterior mesure. At fixed size of the microstructure, we discussed a procedure to account for the modelling error. We also proposed an efficient algorithm to sample from the posterior measure combining numerical homogenization and reduced basis techniques. Several numerical examples illustrating the efficiency of the proposed method and confirming our theoretical findings were also given.

## References

[1] ABDULLE, A. On a priori error analysis of fully discrete heterogeneous multiscale FEM. *Multiscale Model. Simul. 4*, 2 (2005), 447–459.

[2] ABDULLE, A. The finite element heterogeneous multiscale method: a computational strategy for multiscale PDEs. In *Multiple scales problems in biomathematics, mechanics, physics and numerics*, vol. 31 of *GAKUTO Internat. Ser. Math. Sci. Appl.* Gakkōtosho, Tokyo, 2009, pp. 133–181.

[3] ABDULLE, A., AND BAI, Y. Reduced basis finite element heterogeneous multiscale method for high-order discretizations of elliptic homogenization problems. *J. Comput. Phys. 231*, 21 (2012), 7014–7036.

[4] ABDULLE, A., AND DI BLASIO, A. Numerical homogenization and model order reduction for multiscale inverse problems. *Submitted to publication* (2016).

[5] ABDULLE, A., AND NONNENMACHER, A. A short and versatile finite element multiscale code for homogenization problems. *Comput. Methods Appl. Mech. Engrg. 198*, 37–40 (2009), 2839–2859.

[6] ARRIDGE, S. R., KAIPIO, J. P., KOLEHMAINEN, V., SCHWEIGER, M., SOMERSALO, E., TARVAINEN, T., AND VAUHKONEN, M. Approximation errors and model order reduction with an application in optical diffusion tomography. *Inverse Problems 22* (2006), 175–195.

[7] BENSOUSSAN, A., LIONS, J.-L., AND PAPANICOLAOU, G. *Asymptotic analysis for periodic structures.* North-Holland Publishing Co., Amsterdam, 1978.

[8] BOGACHEV, V. I. *Measure Theory.* Springer-Verlag, Berlin, Heidelberg, 2007.

[9] CALVETTI, D., ERNST, O., AND SOMERSALO, E. Dynamic updating of numerical model discrepancy using sequential sampling. *Inverse Problems 30*, 11 (2014).

[10] CAREY, G. F., CHOW, S.-S., AND SEAGER, M. K. Approximate boundary-flux calculations. *Computer Methods in Applied Mechanics and Engineering 50*, 2 (1985), 107–120.

[11] CIORANESCU, D., AND DONATO, P. *An introduction to homogenization*, vol. 17 of *Oxford Lecture Series in Mathematics and its Applications.* Oxford University Press, New York, 1999.

[12] DASHTI, M., AND STUART, A. M. The bayesian approach to inverse problems. In *Handbook of Uncertainty Quantification.* Springer, 2016.

[13] DE GIORGI, E., AND SPAGNOLO, S. Sulla convergenza degli integrali dell'energia per operatori ellittici del secondo ordine. *Boll. Un. Mat. Ital. 4*, 8 (1973), 391–411.

[14] DUNLOP, M. M., AND STUART, A. M. The bayesian formulation of EIT: analysis and algorithms. 1007–1036.

[15] GREPL, M. A., MADAY, Y., NGUYEN, N. C., AND PATERA, A. T. Efficient reduced-basis treatment of nonaffine and nonlinear partial differential equations. *ESAIM: Mathematical Modelling and Numerical Analysis-Modélisation Mathématique et Analyse Numérique 41*, 3 (2007), 575–605.

[16] IGLESIAS, M. A., LU, Y., AND STUART, A. M. A bayesian level set method for geometric inverse problems. 181–217.

[17] MURAT, F., AND TARTAR, L. $H$-convergence. In *Topics in the mathematical modelling of composite materials*, vol. 31 of *Progr. Nonlinear Differential Equations Appl.* Birkhäuser Boston, Boston, MA, 1997, pp. 21–43.

[18] NOLEN, J., PAVLIOTIS, G. A., AND STUART, A. M. Multiscale modelling and inverse problems. In *Numerical Analysis of Multiscale Problems*, I. G. Graham, T. Y. Hou, O. Lakkis, and R. Scheichl, Eds., vol. 83 of *Lecture Notes in Computational Science and Engineering.* Springer-Verlag Berlin Heidelberg, 2012, pp. 1–34.

[19] STUART, A. M. Inverse problems: A bayesian perspective. 451–559.

[20] TARTAR, L. *Estimations des coefficients homogénéisés.* Lectures Notes in Mathematics 704. Springer-Verlag, Berlin, 1977.